\documentclass[12pt]{article}
\usepackage{amsmath, amsfonts, amscd, latexsym, amsthm, amssymb, graphicx}
\newtheorem{theor}{\hspace{1cm}{\sc Theorem}}[section]
\newtheorem{theor0}{\hspace{1cm}{\sc Theorem}}
\newtheorem{utver}[theor]{\hspace{1cm}{\sc Assertion}}
\newtheorem{sledst}[theor]{\hspace{1cm}{\sc Corollary}}
\newtheorem{lemma}[theor]{\hspace{1cm}{\sc Lemma}}
\newtheorem{conjec}[theor]{\hspace{1cm}{\sc Conjecture}}
\theoremstyle{definition}
\newtheorem{defin}[theor]{\hspace{1cm}{\sc Definition}}
\newtheorem{exa}[theor]{\hspace{1cm}{\sc Example}}
\newtheorem*{rem}{\hspace{1cm}{\sc Remark}}

\newcommand{\ind}{\mathop{\rm ind}\nolimits}

\newcommand{\Span}{\mathop{\rm conv}\nolimits}

\newcommand{\Vol}{\mathop{\rm Vol}\nolimits}

\renewcommand{\Re}{\mathop{\rm \rm Re}\nolimits}
\newcommand{\codim}{\mathop{\rm codim}\nolimits}
\newcommand{\conv}{\mathop{\rm conv}\nolimits}

\newcommand{\Int}{\mathop{\rm Int}\nolimits}
\newcommand{\id}{\mathop{\rm id}\nolimits}
\newcommand{\res}{\mathop{\rm res}\nolimits}

\newcommand{\rk}{\mathop{\rm rk}\nolimits}
\newcommand{\Rez}{\mathop{\rm Res}\nolimits}
\newcommand{\pt}{\mathop{\rm pt}\nolimits}

\def\R{\Bbb R}
\def\N{\Bbb N}
\def\Z{\Bbb Z}
\def\Q{\Bbb Q}
\def\C{\Bbb C}
\def\CC{({\Bbb C}\setminus 0)}
\def\T{\Bbb T}

\def\CP{\Bbb C\Bbb P}

\def\K{\Bbb K}

\begin{document}

\begin{center}
{\Large Determinantal singularities and Newton
polyhedra\footnote{ For the most part, this is a significantly
updated exposition of material from \cite{E1}-\cite{E4}.}}
\end{center}

\begin{center}
{A. I. Esterov}\footnote{ Supported in part by
RFBR-JSPS-06-01-91063, RFBR-07-01-00593, INTAS-05-7805.}
\end{center}

There are two well researched tasks, related to Newton polyhedra:
\begin{quote}
A) to describe Newton polyhedra of resultants and discriminants;

B) to study invariants of singularities in terms of their Newton
polyhedra.
\end{quote}
We introduce so called resultantal singularities (Definitions
\ref{torrezcikl} and \ref{defressing}), whose study in terms of
Newton polyhedra unifies the tasks A and B to a certain extent.
In particular, this provides new formulations and proofs for a
number of well known results (see, for example, Theorem
\ref{thsturmf} related to task A and Corollary \ref{oka1f}
related to task B).

As an application, we study basic topological invariants of
determinantal singularities (Subsection \ref{ssdet}) and certain
generalizations of the Poincare-Hopf index (Subsection
\ref{ssind}) in terms of Newton polyhedra. By generalizations of
the Poincare-Hopf index we mean invariants of collections of
(co)vector fields that participate in generalizations of the
classical Poincare-Hopf formula to singular varieties, arbitrary
characteristic numbers, etc.

In order to study resultantal singularities in terms of Newton
polyhedra, we introduce relative versions of the mixed volume
(Subsection \ref{ssmixvol}), Kouchnirenko--Bernstein formula
(Subsection \ref{sskb}), Kouchnirenko--Bernstein--Kho\-van\-skii
formula (Subsection \ref{sskbk}), construct a toric resolution for
a resultantal singularity (Subsection \ref{sstorrez}), and
introduce one more tool, which has no counterpart for complete
intersection singularities (Subsection \ref{sscay}).

Recall that the Kouchnirenko--Bernstein formula \cite{bernst}
expresses the number of solutions of a system of polynomial
equations $f_1(x_1,\ldots,x_n)=\ldots=f_n(x_1,\ldots,x_n)=0$ in
terms of the Newton polyhedra of the polynomials $f_1,\ldots,f_n$,
provided that the coefficients of $f_1,\ldots,f_n$ satisfy a
certain condition of general position (the \textit{Newton
polyhedron} of a polynomial is the convex hull of the exponents
of its monomials).

A local version of this formula is presented below as an
illustrative special case of our results. Let $P$ be the set of
all polyhedra $\Delta$ in the positive orthant
$\R^n_+\subset\R^n$, such that the difference
$\R^n_+\setminus\Delta$ is bounded. $P$ is a semigroup with
respect to the Minkowski addition $\Delta_1+\Delta_2=\{x+y\;|\;
x\in\Delta_1,\, y\in\Delta_2\}$. The symmetric multilinear
function $\Vol:\underbrace{P\times\ldots\times P}_n\to\R$, such
that $\Vol(\Delta,\ldots,\Delta)$ equals the volume of
$\R^n_+\setminus\Delta$ for every $\Delta\in P$, is called the
\textit{mixed volume} of $n$ polyhedra. The minimal polyhedron in
$P$ that contains the exponents of all monomials, participating in
a power series $f$ of $n$ variables, is called the \textit{Newton
polyhedron} of $f$ (provided that such minimal polyhedron exists).
The coefficient of a monomial in the power series $f$ is called a
\textit{leading coefficient}, if the exponent of this monomial is
contained in a bounded face of the Newton polyhedron. For a
linear function $l:\Z^n\to\Z$ and a power series $f$ in $n$
variables, denote the lowest order non-zero $l$-quasihomogeneous
component of $f$ by $f^l$.
\begin{theor0}[\cite{E4}] \label{th1} If $\Delta_1,\ldots,\Delta_n$ are the Newton
polyhedra of complex analytic germs $f_1,\ldots,f_n,\;
f_i:(\C^n,0)\to(\C,0)$, then the topological degree of the map
$(f_1,\ldots,f_n):(\C^n,0)\to(\C^n,0)$ is greater or equal to
$n!\Vol(\Delta_1,\ldots,\Delta_n)$.

If, in addition, the leading coefficients of $f_1,\ldots,f_n$
satisfy a certain condition of general position, then the
topological degree equals $n!\Vol(\Delta_1,\ldots,\Delta_n)$.

The necessary and sufficient condition of general position is
that, for every linear function $l:\Z^n\to\Z$ with positive
coefficients, the polynomial equations $f^l_1=\ldots=f^l_n=0$ have
no common roots in $\CC^n$.
\end{theor0}
The proof is the same as for the classical D.~Bernstein's formula
(although it seems that the right hand side of the local version
was never described in this form before, cf. \cite{oka} and
\cite{biviam}). We generalize it in the following directions.

\textbf{(1)} The topological degree of the map
$(f_1,\ldots,f_n):(\C^n,0)\to(\C^n,0)$ can be regarded as the
Poincare-Hopf index of the complex vector field $(f_1,\ldots,f_n)$
on $\C^n$. We extend Theorem \ref{th1} to certain generalizations
of the Poincare-Hopf index (\cite{smg}, \cite{smgfam},
\cite{smgdet}, \cite{suwa}, etc), see Subsection \ref{ssind}.

\textbf{(2)} The topological degree of the map $(f_1,\ldots,f_n)$
differs by 1 from the Milnor number of the 0-dimensional complete
intersection $f_1=\ldots=f_n=0$. We study Milnor numbers of
arbitrary-dimensional complete intersections and, more generally,
of determinantal singularities, see Subsection \ref{ssdet}.

\textbf{(3)} The topological degree of the map $(f_1,\ldots,f_n)$
can be regarded as the intersection number of the divisors
$\{f_1=0\},\ldots,\{f_n=0\}$. More generally, we may assume that
$f_1,\ldots,f_n$ are sections of line bundles on an arbitrary
toric variety, see Subsection \ref{sskb}.

Subsections \ref{ssmixvol}--\ref{ssdet}, \ref{sstor}--\ref{sskb}
and \ref{ssresvar} introduce necessary notation and recall some
basic facts related to convex geometry (Subsection
\ref{ssmixvol}), mixed volumes (Subsection \ref{ssmixvol}),
resultants (Subsections \ref{ssres} and \ref{ssresvar}),
invariants of singularities (Subsection \ref{sstopinv}), Newton
polyhedra (Subsections \ref{ssdet} and \ref{sstor}), toric
varieties (Section \ref{sstor}) and intersection theory (Section
\ref{sskb}).

This work is based on the author's thesis, I am grateful to
Professor S.~M.~Gusein-Zade for his ideas and guidance. I also
want to thank A.~G.~Khovanskii and S.~P.~Chulkov for many
valuable remarks.
\vspace{-0.5cm}

{\scriptsize \tableofcontents}


\section{Applications}\label{ss1} Before discussing the main results of
the paper, we present some of their applications in this section.
In the first subsection, we introduce mixed volume of pairs of
polytopes. In the second subsection, we present a formula for the
support function of the Newton polytope of the resultant in terms
of mixed volumes of pairs. In the third subsection, we recall
definitions of invariants of singularities that we wish to count
in terms of Newton polyhedra. In the last two subsections, we
count these invariants in terms of mixed volumes of Newton
polyhedra of singularities, provided that principal parts of
singularities are in general position.
\subsection{Relative mixed volume and prisms}\label{ssmixvol}

There are two possible settings in this subsection:

I) A \textit{polyhedron} in $\R^n$ is an intersection of finitely
many closed half-spaces. The \textit{volume form} in $\R^n$ and in
every its subspace is induced by the standard metric in $\R^n$.
The set $S\subset(\R^n)^*$ consists of all unit covectors in the
sense of the standard metric in $\R^n$, and $\K$ stands for $\R$.

II) A \textit{polyhedron} in $\R^n$ is an intersection of
finitely many rational closed half-spaces, whose vertices are
contained in $\Z^n$. The \textit{volume form} in $\R^n$ and every
its rational subspace is chosen in such a way that the minimal
possible volume of a parallelepiped with integer vertices equals
1. The set $S\subset(\R^n)^*$ consists of all primitive covectors
(an integer covector is said to be \textit{primitive} if it is
not equal to another integer covector multiplied by a positive
integer number). $\K$ stands for $\frac{\Z}{n!}$.

Let $\mathcal{M}$ be the set of all convex bounded polyhedra in
$\R^n$. This set is a semigroup with respect to the Minkowski
addition $A+B=\{a+b\; |\; a\in A,\; b\in B\}$. The \textit{mixed
volume} is the unique symmetric multilinear function
$$\Vol:\underbrace{\mathcal{M}\times\ldots\times
\mathcal{M}}_n\to\K,$$ such that $\Vol(A,\ldots,A)$ equals the
volume of $A$ for every polyhedron $A\in \mathcal{M}$.

We introduce a "relative" version of the mixed volume. Let
$N\subset\R^n$ be a convex polyhedron (not necessary bounded).
Its \textit{support function} $N(\cdot)$ is defined as
$$N(\gamma)=\inf\limits_{x\in N} \gamma(x)$$ for every covector
$\gamma\in (\R^n)^*$. The set $$N^{\gamma}=\{x\in N\; |\;
\gamma(x)=N(\gamma)\}$$ is called the \textit{support face} of the
polyhedron $N$ with respect to the covector $\gamma\in (\R^n)^*$.
The set $\{\gamma\; |\; N(\gamma)>-\infty\}\subset (\R^n)^*$ is
called the \textit{support cone} of $N$.

Consider the set $\mathcal{M}_{\Gamma}$ of all ordered pairs of
polyhedra $(A,B)$ with a given support cone
$\Gamma\subset(R^n)^*$, such that the symmetric difference
$A\vartriangle B$ is bounded. $\mathcal{M}_{\Gamma}$ is a
semigroup with respect to Minkowski addition of pairs
$(A,B)+(C,D)=(A+C,B+D)$. Example of Minkowski addition:

\begin{center}
\noindent\includegraphics[width=12cm]{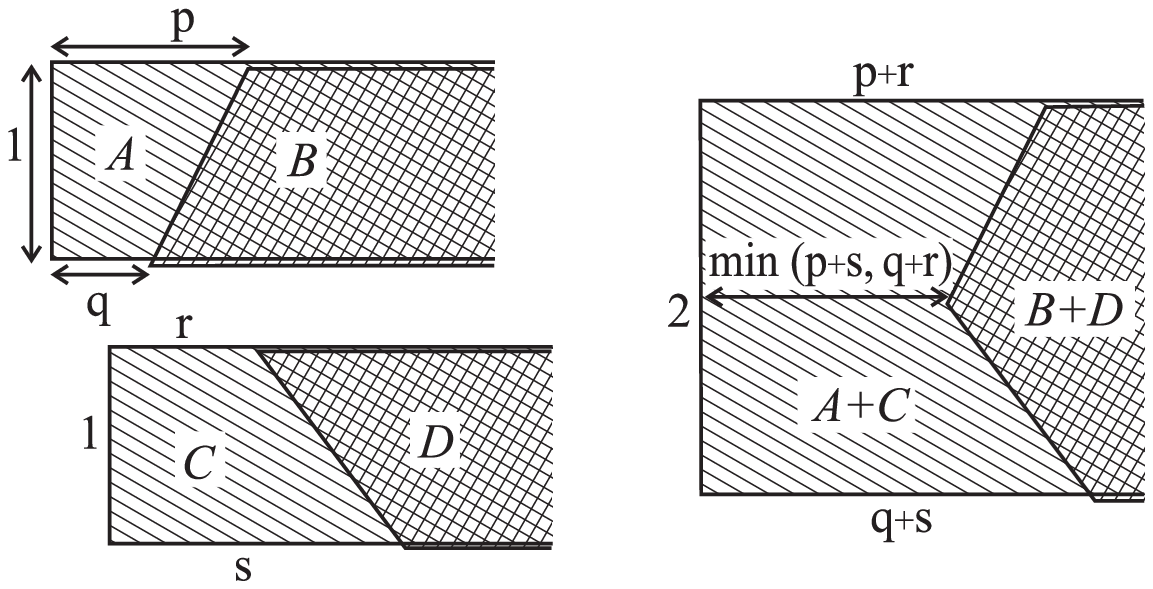} \end{center}

\begin{defin}[\cite{E2}, \cite{E3}] \label{defrelmixvol}
\textit{The volume} $V(A,B)$ of the pair of polyhedra
$(A,B)\in\mathcal{M}_{\Gamma}$ is the difference of the volumes
of the sets $A\setminus B$ and $B\setminus A$. \textit{The mixed
volume} of pairs of polyhedra with the support cone
$\Gamma\subset(\R^n)^*$ is the symmetric multilinear function
$\Vol_{\Gamma}:\underbrace{\mathcal{M}_{\Gamma}\times\ldots\times
\mathcal{M}_{\Gamma}}_n\to\K$, such that
$\Vol_{\Gamma}\bigl((A,B),\ldots,(A,B)\bigr)=V(A,B)$ for every
pair $(A,B)\in \mathcal{M}_{\Gamma}$.
\end{defin}
For example, if $\Gamma=(\R^n)^*$, then $\mathcal{M}_{\Gamma}$ is
the set of pairs of convex bounded polyhedra, and the mixed
volume of pairs
$\Vol_{\Gamma}\Bigl((A_1,B_1),\ldots,(A_n,B_n)\Bigr)$ equals the
difference of the classical mixed volumes of the collections
$A_1,\ldots,A_n$ and $B_1,\ldots,B_n$. The following theorem
provides existence, uniqueness and basic formulas for computation
of the mixed volume of pairs.

\begin{utver}[\cite{E3}] \label{voldef}
1) If a function
$F:\underbrace{\mathcal{M}_{\Gamma}\times\ldots\times
\mathcal{M}_{\Gamma}}_n\to\K$ is symmetric, multilinear, and
$F\bigl((A,B),\ldots,(A,B)\bigr)=V(A,B)$ for every pair $(A,B)\in
\mathcal{M}_{\Gamma}$, then
$$F\bigl((A_1,B_1),\ldots,(A_n,B_n)\bigr)=\frac{1}{n!}\sum\limits_{I\in\{1,\ldots,n\}}
(-1)^{n-|I|} V\bigl(\sum\limits_{i\in I} (A_i,B_i)\bigr).$$
\newline 2) The function
$\Vol_{\Gamma}:\underbrace{\mathcal{M}_{\Gamma}\times\ldots\times
\mathcal{M}_{\Gamma}}_n\to\K$, defined by the equality
$$\Vol_{\Gamma}\bigl((A_1,B_1),\ldots,(A_n,B_n)\bigr)=\frac{1}{n!\cdot n}
\sum_{\sigma\in\mathcal{S}^n} \sum_{k=1}^n
\sum\limits_{\gamma\in\Int\Gamma\cap S}
\bigl(B_{\sigma(k)}(\gamma)-A_{\sigma(k)}(\gamma)\bigr)\times$$
$$\times
\Vol(A_{\sigma(1)}^{\gamma},\ldots,A_{\sigma(k-1)}^{\gamma},
B_{\sigma(k+1)}^{\gamma},\ldots,B_{\sigma(n)}^{\gamma}),$$ is
symmetric, multilinear, and
$\Vol_{\Gamma}\bigl((A,B),\ldots,(A,B)\bigr)=V(A,B)$ for every
pair $(A,B)\in \mathcal{M}_{\Gamma}$ (here $\Vol$ is the classical
mixed volume of bounded polyhedra, $\Int\Gamma$ stands for the
interior of the cone $\Gamma$, and the sum in the right hand side
has only finitely many non-zero terms, because covectors
$\gamma$, corresponding to the non-zero terms, are the external
normal covectors of the bounded codimension 1 faces of the
polyhedron $\sum_{i=1}^n A_i+B_i$).
\end{utver}
\begin{exa} \label{exavol} The mixed volume of pairs $(A,B)$ and $(C,D)$ on the picture above equals $\min(p+s,\, q+r)/2$. \end{exa}
\textsc{Proof} of part 1: substitute every term of the form
$V(A,B)$ in the right hand side of the desired equality by
$F\bigl((A,B),\ldots,(A,$ $B)\bigr)$, open the parentheses by
linearity of $F$, and cancel like terms by the symmetric property
of $F$. Part 2 follows from the linearity and the symmetric
property of the classical mixed volume and the following fact (see
\cite{khovdan}): the volume of the "trapezoid" $A\times\{0\}\cup
B\times\{h\}\subset\R^{n-1}\times\R$ equals
$$\frac{h}{n}\sum_{k=0}^{n-1}
\Vol(\underbrace{A,\ldots,A}_k,\underbrace{B,\ldots,B}_{n-k-1})$$
for bounded polyhedra $A$ and $B$ in $\R^{n-1}.\;\Box$

The symmetrizing in the right hand side of Part 2 (i.e. summation
over all $\sigma\in\mathcal{S}^n$) turns out to be redundant,
which significantly simplifies computations. In particular, the
mixed volume of pairs of polyhedra with integer vertices is
contained in $\frac{\Z}{n!}$.

\begin{utver}[\cite{E3}] \label{volsymm} We have
$\Vol_{\Gamma}\bigl((A_1,B_1),\ldots,(A_n,B_n)\bigr)=$
$$=\frac{1}{n}
\sum_{k=1}^n \sum\limits_{\gamma\in\Int\Gamma\cap S}
\bigl(B_k(\gamma)-A_k(\gamma)\bigr) \Vol(A_1^{\gamma},\ldots,
A_{k-1}^{\gamma},B_{k+1}^{\gamma},\ldots,B_n^{\gamma}).$$
\end{utver}
The proof is based on relations between mixed volumes and
algebraic geometry, and will be given together with the proof of
Theorem \ref{relbernst} in Subsection \ref{ss2proof}. One can
easily verify that the right hand side of the Assertion
\ref{volsymm} is not symmetric for pairs of polyhedra with
unbounded symmetric difference, and its symmetrization is not
contained in $\frac{\Z}{n!}$.

One more formula for the mixed volume of pairs $(A_1,B_1),$
$\ldots,$ $(A_n,B_n)$ will also be proved together with Theorem
\ref{relbernst} in Subsection \ref{ss2proof}:
\begin{utver}[\cite{E3}] \label{volexpr2}
If $\tilde A_i$ and $\tilde B_i$ are bounded polyhedra such that
$A_i\setminus B_i = \tilde A_i\setminus \tilde B_i$ and
$B_i\setminus A_i = \tilde B_i\setminus \tilde A_i$, then
$\Vol_{\Gamma}\bigl((A_1,B_1),\ldots,(A_n,B_n)\bigr)=\Vol(\tilde
A_1,\ldots, \tilde A_n) - \Vol(\tilde B_1,\ldots, \tilde B_n)$.
\end{utver}

We also need the following well known formula:
\begin{lemma} \label{comput1}
Let $A_1,\ldots,A_p$ be bounded polyhedra in a
$(p+q)$-dimensional space $P$, and $B_1,\ldots,B_q$ be bounded
polyhedra in its $q$-dimensional subspace $Q$. Then, denoting the
projection $P\to P/Q$ by $\pi$, we have
$$(p+q)!\Vol(A_1,\ldots,A_p,B_1,\ldots,B_q)=p!q!\cdot\Vol(\pi
A_1,\ldots,\pi A_p)\cdot\Vol(B_1,\ldots,B_q).$$
\end{lemma}

When studying determinantal singularities in
terms of Newton polyhedra, we often deal with mixed volumes of
prisms over Newton polyhedra:
\begin{defin} \label{prismdef} \textit{The prism} $\Delta_1*\ldots*\Delta_n$ \textit{over polyhedra}
$\Delta_1,\ldots,\Delta_n\subset\R^m$ is the convex hull of the
union
$$\bigcup_i \{b_i\}\times\Delta_i\subset\R^{n-1}\oplus\R^m,$$
where $b_1,\ldots,b_n$ are the vertices of the standard simplex in
$\R^{n-1}$. \textit{The prism}
$(\Gamma_1,\Delta_1)*\ldots*(\Gamma_n,\Delta_n)$ \textit{over
pairs of polyhedra}
$(\Gamma_1,\Delta_1),\ldots,(\Gamma_n,\Delta_n)$ in $\R^m$ is the
pair $(\Gamma_1*\ldots*\Gamma_n,\Delta_1*\ldots*\Delta_n)$.
\end{defin} The following formula simplifies the computation of
the mixed volume of integer prisms. For a bounded set
$\Delta\subset\R^m$, denote the number of integer lattice points
in $\Delta$ by $I(\Delta)$. If the symmetric difference of
(closed) integer polyhedra $\Gamma$ and $\Delta$ in $\R^m$ is
bounded, denote the difference
$I(\Gamma\setminus\Delta)-I(\Delta\setminus\Gamma)$ by
$I(\Gamma,\Delta)$. Denote the convex hull of the union of
polyhedra $\Delta_i\subset\R^m$ by $\bigvee_i\Delta_i$. For pairs
of polyhedra $(\Gamma_i,\Delta_i)$ in $\R^m$, denote the pair
$(\bigvee_i\Gamma_i,\bigvee_i\Delta_i)$ by
$\bigvee_i(\Gamma_i,\Delta_i)$. \begin{theor} \label{thint} If
bounded integer polyhedra or pairs of integer polyhedra $B_{i,j}$,
$i=1,\ldots,n$, $j=1,\ldots,k$ in $\R^m$ have the same support
cone, and $m=k-n+1$, then the mixed volume of the prisms
$B_{1,j}*\ldots*B_{n,j}$, $j=1,\ldots,k$, equals
$$\frac{1}{k!}\sum_{J\subset\{1,\ldots,k\}\atop b_1+\ldots+b_{n}=|J|}(-1)^{k-|J|}I\Bigl(\bigvee_{J_1\sqcup\ldots\sqcup J_n=J\atop |J_1|=b_1,\ldots,|J_n|=b_n}
\sum_{i=1,\ldots,n\atop j\in J_i} B_{i,j}\Bigr).$$\end{theor} The
proof is given in Subsection \ref{ss4}. Note that some of
$B_{i,j}$ may be empty.

\subsection{Newton polyhedra of resultants.}\label{ssres}
We denote the monomial $t_1^{a_1}\ldots t_N^{a_N}$ by $t^a$. For
a subset $\Sigma\subset\R^N$, we denote the set of all Laurent
polynomials of the form $\sum\limits_{a\in\Sigma\cap\Z^N} c_a
t^a,\, c_a\in\C$, by $\C[\Sigma]$. We regard them as functions on
the complex torus $(\C\setminus\{0\})^N$.
\begin{defin}
\label{defgkzresultant} For finite sets $\Sigma_i\subset \Z^N,
i=0,\ldots,N$, the resultant $\Rez$ is defined as the equation of
the closure of the set
$$\bigl\{ (g_0,\ldots,g_N)\; |\; g_i\in\C[\Sigma_i],\;
g_0(t)=\ldots=g_N(t)=0 \mbox{ for some } t\in(\C\setminus\{0\})^N
\bigr\},$$ provided that this set is a hypersurface in
$\C[\Sigma_0]\oplus\ldots\oplus\C[\Sigma_N]$. Otherwise, we set
$\Rez=1$ by definition.
\end{defin}
\begin{exa} If $\Sigma_0$ and $\Sigma_1$ are segments in $\Z^1$,
then $\Rez$ is the classical resultant of univariate polynomials
$g_0$ and $g_1$.

If $\Sigma_0=\ldots=\Sigma_N$ is the set of vertices of the
standard $N$-dimensional simplex, then $\Rez$ is the determinant
of the coefficient matrix of an affine linear map
$(g_0,\ldots,g_N)$.
\end{exa}
Without loss in generality, we can assume that the collection
$\Sigma_0,\ldots,\Sigma_N\subset\Z^N$ is \textit{essential}, i.e.
the dimension of the convex hull of $\sum_{j\in J}\Sigma_j$ is not
smaller than $|J|$ for every $J\varsubsetneq\{0,\ldots,N\}$, and
equals $N$ for $J=\{0,\ldots,N\}$ (see \cite{sturmf} for
details). Under this assumption we have the following formula for
the support function of its Newton polyhedron $\Delta_{\Rez}$.

The resultant $\Rez$ is a polynomial in the coefficients
$c_{a,i}$ of indeterminate polynomials $g_i=\sum_{a\in\Sigma_i}
c_{a,i} t^a\in\C[\Sigma_i]$, thus $\Delta_{\Rez}$ is contained in
$\R\otimes\Lambda$, where $\Lambda$ is the lattice of monomials
$\prod c_{a,i}^{\lambda_{a,i}}$. The exponents $\lambda_{a,i}$
can be regarded as linear functions on $\Lambda$ and form a basis
of the dual lattice $\Lambda^*$. For a linear function
$\gamma=\sum\gamma_{a,i}\lambda_{a,i}\in\Lambda^*$ with
coordinates $\gamma_{a,i}$ in this basis, denote the ray
$\{(t,0),\ |\, t\leqslant 0\}\subset\R\oplus\R^N$ by $l$ and the
convex hull $\conv\{(\gamma_{a,i},a)\, |\, a\in
A_i\}+l\subset\R\oplus\R^N$ by $A_{i,\gamma}$.

\begin{theor}[\cite{E2}, \cite{E3}] \label{thsturmf}
The value of the support function $\Delta_{\Rez}(\gamma)$ equals
the mixed volume of pairs
$(N+1)!\cdot\Vol_l\Bigl((A_{0,0},A_{0,\gamma}),\ldots,(A_{N,0},A_{N,\gamma})\Bigr)$.
\end{theor}
The proof is given in Subsection \ref{ss3proof}.
\begin{rem}
When discussing Newton polyhedra of polynomials rather than
Newton polyhedra of analytic germs, the value of the support
function $\Delta(\cdot)$ at a covector $\gamma$ is often defined
as the maximal value of $\gamma$ on the polyhedron $\Delta$,
rather than the minimal one. In this notation, the formula is as
follows: $\Delta_{\Rez}(-\gamma)$ equals
$(N+1)!\cdot\Vol_l\Bigl((A_{0,\gamma},A_{0,0}),\ldots,(A_{N,\gamma},A_{N,0})\Bigr)$.
\end{rem}

This theorem gives a new proof to a number of well known facts
about Newton polyhedra of resultants and discriminants, including
the description of the vertices of the Newton polyhedron
$\Delta_{\Rez}$ (\cite{sturmf}) and the formula for the support
function of the Newton polyhedron of the $A$-determinant
(\cite{gkz}).

\subsection{Topological invariants of singularities}\label{sstopinv}

We recall the definition of singularities and their topological
invariants that we wish to study in terms of Newton polyhedra.

\textbf{Milnor fiber and radial index.}
\begin{defin} \label{defmult} The \textit{multiplicity} of a positive-dimensional complex
analytic germ $V\subset\C^m$ is defined as its intersection
number with a generic vector subspace of complementary dimension
in $\C^m$.
\end{defin}

Suppose that the germ of a complex analytic set $V\subset\C^m$ is
smooth outside the origin. Let $f:(\C^m,0)\to(\C,0)$ be a germ of
a complex analytic function, such that the restriction $f|_{V}$
has no singular points in a punctured ball $B$, centered at the
origin.
\begin{defin}
For a small $\delta\ne 0$, the manifold $V\cap B\cap
f^{(-1)}(\delta)$ is called the \textit{Milnor fiber} of the germ
$f|_V$.
\end{defin}

Suppose that the germ of a real analytic set $V\subset\R^m$ is
smooth outside the origin. Suppose that $\omega$ is a germ of a
real continuous 1-form in $\R^m$ near the origin, and the
restriction $\omega|_V$ has no zeroes in a punctured neighborhood
$U$ of the origin. Let $\tilde \omega$ be a 1-form on
$V\setminus\{0\}$ such that
\newline 1) $\tilde\omega$ has isolated zeroes $p_1,\ldots, p_N$
in $U$, \newline 2) $\tilde\omega=\omega|_V$ outside $U$, and
\newline 3) $\tilde\omega(x)=d(\|x\|^2)|_V$ near the origin.
\begin{defin}(\cite{smg}) The \textit{radial index} of $\omega|_V$ is defined as $1+\sum_j \ind_{p_j}\tilde\omega$,
where $\ind_{p_j}$ is the Poincare-Hopf index at $p_j$.
\end{defin}
If $V$ is smooth, then the radial index of $\omega|_V$ equals the
Poincare-Hopf index $\ind_0 \omega|_V$. If $V\subset\C^m$ and
$f:\C^m\to\C$ are complex analytic, then the sum of the radial
index of the 1-form $d\Re f|_V$ and the Euler
characteristic of a Milnor fiber of $f|_V$ equals 1, see \cite{smgfam}.

\textbf{Determinantal singularities.} Suppose that $A=(a_{i,j}) :
\C^n\to\C^{I\times k}$ is a germ of a matrix with holomorphic
entries near the origin, and $I\leqslant k$ (we denote the space
of all $(I\times k)$-matrices by $\C^{I\times k}$).

\begin{defin} The set $[A]=\{ x\, |\, \rk A(x)<I \}$ is called the
\textit{$(I\times k)$-determinantal set, defined by $A$,}
provided that its dimension equals $n-k+I-1$ (i. e. the minimal
possible one).
\end{defin}

\begin{exa} An $(1\times k)$-determinantal set is a complete intersection of codiminsion $k$. \end{exa}

The following topological invariants of determinantal
singularities and complete intersections can be regarded as
generalizations of the Poincare-Hopf index. We denote the set of
all $I\times k$ degenerate complex matrices by
$\mathcal{D}^{I\times k}\subset\C^{I\times k}$.
\begin{defin}
\label{defegsect} Suppose that the entries of an $I_j\times k_j$
matrix $W_j$ are germs of continuous functions on $(\C^n,0)$,
where $j$ ranges from $1$ to $J$, and $n=\sum_j 1+|k_j-I_j|$. If
the origin is the only point of $\C^n$, where the matrices
$W_1,\ldots,W_J$ are all degenerate, then the intersection number
of the graph of the mapping
$W=(W_1,\ldots,W_J):\C^n\to\C^{I_1\times
k_1}\oplus\ldots\oplus\C^{I_J\times k_J}$ with the product
$\C^n\times\mathcal{D}^{I_1\times
k_1}\times\ldots\times\mathcal{D}^{I_J\times k_J}$ is called
\textit{the multiplicity of the collection} $(W_1,\ldots,W_J)$.
\end{defin}
In other words, if generic small perturbations of the matrices
$W_1,\ldots,W_J$ degenerate at finitely many points near the
origin, then the number of these points (counted with appropriate
signs) equals the multiplicity of $(W_1,\ldots,W_J)$. Since
determinantal singularities are Cohen-Macaulay, the multiplicity
of a holomorphic collection $(W_1,\ldots,W_J)$ equals
$$\dim_{\C}\mathcal{O}_{(\C^n,0)}/\langle\mbox{maximal minors
of}\; W_1,\ldots,W_J\rangle$$ \noindent(it is also referred to as
the Buchsbaum-Rim multiplicity in this case). Note that the
multiplicity of the collection $(W_1,\ldots,W_J)$ is not equal to
the intersection number of the image $M=W(\C^n)$ and the product
$\mathcal{D}^{I_1\times
k_1}\times\ldots\times\mathcal{D}^{I_J\times k_J}$ in general ---
their ratio is equal to the topological degree of the germ of the
mapping $W:\C^n\to M$.

Definition \ref{defegsect} can be regarded as a generalization of the 
Poincare-Hopf index due to the following observation. Let
$v_{i,j}$ be a smooth section of a vector bundle $\mathcal{I}_j$
of rank $k_j$ on a smooth $n$-dimensional complex manifold $M$ for
$i=1,\ldots,I_j,\; I_j\leqslant k_j,\; j=1,\ldots,J$, and
$n=\sum_j 1+k_j-I_j$. Suppose that, for every point $x\in M$,
except for a finite set $X\subset M$, there exists $j$ such that
the vectors $v_{1,j}(x),\ldots,v_{I_j,j}(x)$ are linearly
independent. Choosing a local basis $s_{1,j},\ldots,s_{k_j,j}$ of
the bundle $\mathcal{I}_j$ near a point $x\in X$, one can
represent $v_i$ as a linear combination
$v_i=w_{i,1,j}s_1+\ldots+w_{i,k_j,j}s_{k_j}$, where
$w_{\cdot,\cdot,j}$ are the entries of a smooth $(I_j\times
k_j)$-matrix $W_j:M\to\C^{I_j\times k_j}$, defined near $x$.
Denote the multiplicity of the collection $(W_1,\ldots,W_J)$ by
$m_x$. Then the Chern number
$c_{1+k_1-I_1}(\mathcal{I}_1)\smile\ldots\smile
c_{1+k_J-I_J}(\mathcal{I}_J)\cdot[M]$ is equal to the sum of the
multiplicities $m_x$ over all points $x\in X$, which is the
classical Poincare-Hopf formula for $J=I_1=1$ (see, for example,
\cite{gh}).

The special case of Definition \ref{defegsect} for $I_1=1,\; J=2$,
is called \textit{the Suwa residue} of the collection of sections
$(w^{1,2},\ldots,w^{I_2,2})^T$ of a $k_2$-dimensional vector
bundle on a germ of the complete intersection $w^{1,1}=0$ (see
\cite{suwa}). Another special case is the Gusein-Zade--Ebeling
index of a 1-form on a complete intersection.

\begin{defin}[\cite{smg}, \cite{smgfam}] \label{defegform}
Let $f_1,\ldots,f_k$ be germs of holomorphic functions on
$(\C^n,0)$. Suppose that $f_1=\ldots=f_k=0$ is an isolated
singularity of a complete intersection, i.e. the 1-forms
$df_1,\ldots,df_k$ are linearly independent at every point of the
set $\{f_1=\ldots=f_k=0\}\setminus\{0\}$. Let $\omega$ be a germ
of a smooth 1-form on $(\C^n,0)$, whose restriction to
$\{f_1=\ldots=f_k=0\}\setminus\{0\}$ has no zeros near the
origin. \textit{The Gusein-Zade--Ebeling index} of $\omega$ on the
complete intersection $f_1=\ldots=f_k=0$ is defined as the
multiplicity of the collection of the two matrices
$$\begin{pmatrix}\omega \\ df_1 \\ \vdots \\ df_k\end{pmatrix},
(f_1, \ldots, f_k).$$
\end{defin}

This can be regarded as a generalization of the Poincare-Hopf
index, because the sum of the indices of a 1-form on a variety
with isolated complete intersection singularities equals the
Euler characteristic of the smoothing of the variety.

\subsection{Determinantal singularities and Newton polyhedra}\label{ssdet}

We study the aforementioned invariants of determinantal
singularities in terms of Newton polyhedra.

\textbf{Newton polyhedra.} The monomial
$x_1^{a_1},\ldots,x_n^{a_n}$ is denoted by $x^a$, where
$a=(a_1,\ldots,a_n)\in\Z^n$. The positive orthant of $\R^n$ is
denoted by $\R^n_+$.

\begin{defin}[\cite{E1}]
\textit{The Newton polyhedron} $\Delta_f$ of a germ of a function
$f:(\C^n,0)\to (\C,0),\; f(x)=\sum\limits_{a\in A} c_a x^a$,
where $c_a\ne 0$ for all $a\in A$, is defined as the convex hull
of the Minkowski sum $A+\R^n_+$. The coefficient $c_a$ is said to
be a \textit{leading coefficient} of the power series $f$, if $a$
is contained in a bounded face of the Newton polyhedron
$\Delta_f$.

\textit{The Newton polyhedron} $\Delta_{\omega}$ of a germ of a
1-form $\omega=\sum_i\omega_idx_i$ on $(\C^n,0)$ is defined as the
convex hull of the union of the Newton polyhedra
$\Delta_{\omega_ix_i},\; i=1,\ldots,n$. The coefficient of the
monomial $x^adx_i/x_i$ in the power series expansion of $\omega$
is said to be leading, if $a$ is contained in a bounded face of
$\Delta_{\omega}$.
\end{defin}
Note that $\Delta_{df}=\Delta_f$ for every function $f$. Recall
that, for a collection of positive weights $l=(l_1,\ldots,l_n)$,
assigned to the variables $x_1,\ldots,x_n$ and their differentials
$dx_1,\ldots,dx_n$, the lowest order non-zero $l$-quasihomogeneous
component of the power series $f=\sum\limits_{a\in \Z^n} c_a x^a$
and $\omega=\sum\limits_{a\in \Z^n} c_{a,i} x^adx_i$ is denoted by
$f^l$ and $\omega^l$ respectively (in the latter case, the weight
of the differential $dx_i$ is equal to the weight of $x_i$ by
convention).

\textbf{Determinantal singularities and Newton polyhedra.}
Consider a germ of a holomorphic matrix
$A=(a_{i,j}):\C^n\to\C^{I\times k},\; I<k$. Suppose for
simplicity that the Newton polyhedron of the entry $a_{i,j}$ does
not depend on $i$, being equal to a certain polyhedron
$\Delta_j$ for $i=1,\ldots,I$ and $j=1,\ldots,k$. (we refer this
assumption to as the \textit{unmixedness assumption}).

\begin{defin} The leading coefficients of $A$ are said to be \textit{in
general position}, if, for every collection $l$ of positive
weights and every subset $\mathcal{I}\subset\{1,\ldots,I\}$, the
set of all points $x\in\CC^n$, such that the matrix
$\bigl(a^l_{i,j}(x)\bigr)_{i\in\mathcal{I}\quad\;\;\;\,\atop
j\in\{1,\ldots,k\}}$ is degenerate, has the maximal possible
codimension $k-|\mathcal{I}|+1$.
\end{defin}

For a polyhedron $\Delta\subset\R^n_+$, denote the pair
$(\R^n_+,\Delta)$ by $\widetilde\Delta$. Denote the pair
$(\R^n_+,\, \R^n_+\setminus$standard $n$-dimensional simplex$)$
by $L$.

\begin{theor}[ \cite{E4}] \label{matrmult} Denote the Newton polyhedron of the entry
$a_{i,j}$ of a holomorphic matrix $A=(a_{i,j}):\C^n\to\C^{I\times
k},\; I<k$ by $\Delta_j$. Suppose that the differences
$\R^n_+\setminus\Delta_j,\; j=1,\ldots,k$, are bounded and the
leading coefficients of $A$ are in general position. Then $A$
defines the germ of a $(I\times k)$-determinantal set $[A]$, whose
multiplicity equals
$$\sum_{0<j_0<\ldots<j_{k-I}\leqslant k} n!\Vol_{\R^n_+}(\widetilde\Delta_{j_{\,0}},\ldots,\widetilde\Delta_{j_{\,k-I}},\underbrace{L,\ldots,L}_{n-k+I-1}).$$
If the leading coefficients of $A$ are not in general position,
then the multiplicity is greater than expected, or $A$ is not a
determinantal singularity.
\end{theor}
If $\dim[A]>0$, then the
multiplicity should be understood in the sense of Definition
\ref{defmult}, otherwise in the sense of Definition
\ref{defegsect}. In both cases, this theorem follows from Theorem
\ref{egvol}, which also allows to drop the unmixedness
assumption. A similar formula for the Buchsbaum-Rim multiplicity
of the matrix $A$ is given in \cite{bivia} for $n\leqslant k-I+1$.

\begin{defin} The leading coefficients of $A$ are said to be \textit{in
strong general position}, if, for every collection $l$ of positive
weights, the polynomial matrix $(a^l_{i,j})$ defines a nonsingular
determinantal set in $\CC^n$.

In this case, the leading coefficients of a germ of a 1-form
$\omega$ are said to be \textit{in general position with respect
to} $A$, if, for every collection $l$ of positive weights, the
restriction of $\omega^l$ to the determinantal set, defined by
the matrix $(a^l_{i,j})$ in $\CC^n$, has no zeros.
\end{defin}

For a set $\mathcal{J}\subset\{1,\ldots,n\}$, let
$\R^\mathcal{J}\subset\R^n$ be a coordinate plane given by the
equations $x_i=0,\, i\notin \mathcal{J}$. For a polyhedron
$\Delta\subset\R^n_+$, denote the pair
$(\R^\mathcal{J}\cap\R^n_+,\, \R^\mathcal{J}\cap\Delta)$ by
$\widetilde\Delta^\mathcal{J}$.
\begin{theor}[\cite{E4}] \label{matrzeta} Denote the Newton polyhedron of the entry
$a_{i,j}$ of a holomorphic matrix $A=(a_{i,j}):\C^n\to\C^{I\times
k},\; I<k$ by $\Delta_j$. Suppose that the differences
$\R^n_+\setminus\Delta_j,\; j=1,\ldots,k$ are bounded, and
$n\leqslant 2(k-I+2)$.
\newline 1) If the leading coefficients of $A$ are in strong general position, then the
determinantal set $[A]$ is smooth outside the origin.
\newline 2) If the Newton polyhedron $\Delta_0\subset\R^n_+$ of a germ $f:\C^n\to\C$ intersects all coordinate axes,
and the leading coefficients of $df$ are in general position
w.r.t. $A$, then the Euler characteristic of a Milnor fiber of
$f|_{[A]}$ equals
$$\chi(\Delta_0,\ldots,\Delta_k)=\sum_{a_0\in\N,\;
\mathcal{J}\subset\{1,\ldots,n\},\atop\{j_1,\ldots,j_q\}\subset\{1,\ldots,k\}}(-1)^{|\mathcal{J}|+k-I}C_{|\mathcal{J}|+q-a_0-2}^{q-k+I-1}\times$$
$$\times\sum_{a_{j_1}\in\N,\, \ldots,\, a_{j_q}\in\N,\atop
a_{j_1}+\ldots+a_{j_q}=|\mathcal{J}|-a_0}|\mathcal{J}|!\Vol_{\R^n_+}(\underbrace{\widetilde\Delta_0^\mathcal{J},\ldots,\widetilde\Delta_0^\mathcal{J}}_{a_0},
\underbrace{\widetilde\Delta_{j_1}^\mathcal{J},\ldots,\widetilde\Delta_{j_1}^\mathcal{J}}_{a_{j_1}},\ldots,
\underbrace{\widetilde\Delta^\mathcal{J}_{j_q},\ldots,\widetilde\Delta^\mathcal{J}_{j_q}}_{a_{j_q}}).$$
\newline 3) If the Newton polyhedron $\Delta_0\subset\R^n_+$ of a germ of a 1-form $\omega$ on $(\C^n,0)$ intersects all coordinate axes,
and the leading coefficients of $\omega$ are in general position
w.r.t. $A$, then the radial index of $\omega|_{[A]}$ makes sense
and equals $1-\chi(\Delta_0,\ldots,\Delta_k)$.
\end{theor}
In the formula above, $C_n^k=0$ for $k\notin\{0,\ldots,n\}$ by
convention. Thus, all the terms in this sum vanish, except for
those with $|\mathcal{J}|-a_0\geqslant q>k-I$.

The proof is given in Subsection \ref{ss3proof}. Actually, the
same argument allows to drop the unmixedness assumption, and to
compute the $\zeta$-function of monodromy for the restriction
$f|_{[A]}$ (or, more generally, for the restriction of $f$ to an
isolated resultantal singularity), see \cite{E4} for details.
Parts 2 and 3 together provide a formula for the Poincare-Hopf
index of a 1-form on a smoothable determinantal singularity (see
\cite{smgdet}).

The assumption of general position, that we impose on the leading
coefficients of the 1-form, implies that the Newton polyhedra of
its coordinates are almost equal to each other (see the definition
of interlaced polyhedra in \cite{E1}). When computing the radial
index of a 1-form on a complete intersection, one can drop this
assumption as follows: the radial index of a 1-form on a complete
intersection differs from its Ebeling--Gusein-Zade index by the
Milnor number of the complete intersection (see \cite{smg}); the
latter two invariants can be computed in terms of Newton
polyhedra by Theorem \ref{finoka1} and Corollary \ref{oka1f}
respectively. It would be interesting to relax the assumption of
general position in the same way for an arbitrary determinantal
singularity.

\subsection{(Co)vector fields and Newton polyhedra}\label{ssind}

\begin{defin} In the notation of Definition
\ref{defegsect}, the leading coefficients of the collection
$(W_1,\ldots,W_J),\; W_j=(w_k^{i,j})$, are said to be \textit{in
general position}, if, for every collection $l$ of positive
weights, assigned to the variables $x_1,\ldots,x_n$ and
$t_{i,j},\; i=1,\ldots,I_j,\; j=1,\ldots, J$, the lowest order
non-zero $l$-homogeneous components of the linear combinations
$$\sum_{i=1}^{I_j}w^{i,j}_k(x_1,\ldots,x_n)t_{i,j},\; k=1,\ldots,k_j,\;
j=1,\ldots,J,$$ which are polynomials of variables $x_{\cdot}$
and $t_{\cdot,\cdot}$, have no common zeros outside the coordinate
planes.
\end{defin}

\begin{theor}[\cite{E3}]
\label{egvol} Adopting notation of Definition \ref{defegsect},
denote the Newton polyhedron of the $(i,k)$-entry of the matrix
$W_j$ by $\Delta^{i,j}_k$, and suppose that the difference
$\R^n_+\setminus\Delta^{i,j}_k$ is bounded. If the leading
coefficients of the entries of the matrices $W_1,\ldots,W_J$ are
in general position, then the multiplicity of this collection of
matrices makes sense and equals $(k_1+\ldots+k_J)!$ times the
mixed volume of images of prisms
$$(\R^n_+,\Delta^{1,j}_k)*\ldots*(\R^n_+,\Delta^{I_j,j}_k)\subset\R^{I_j-1}\oplus\R^n$$
under the natural inclusions
$$\R^{I_j-1}\oplus\R^n\hookrightarrow\R^{I_1-1}\oplus\ldots\oplus\R^{I_J-1}\oplus\R^n,$$
where $k=1,\ldots,k_j$ and $j=1,\ldots,J$. If the leading
coefficients are not in general position, then the multiplicity
is greater than the expected value or is not defined.
\end{theor}

The proof is given in Subsection \ref{ss3proof}.

\begin{exa} \label{exadet} Consider a holomorphic $(2\times 3)$-matrix valued function $A$ of
two variables $x$ and $y$, such that its rows are equal to
\begin{center} $px^2+qy^3+$(higher order terms) and
$rx^5+sy^4+$(higher order terms),\end{center} where $p,q,r,s$ are
constant row vectors. Then the multiplicity of $A$ at the origin
is at least $34$, which is $3!$ times the volume of the
non-convex polyhedron $\Sigma$ on the picture below, and the
equality is attained if and only if $\det\left(\substack{p\\ q\\
s}\right)\ne 0$ and $\det\left(\substack{p\\ r\\ s}\right)\ne 0$.
The two polygons on the picture are the Newton polygons of the
rows of the matrix $A$.
\begin{center} \noindent\includegraphics[width=14cm]{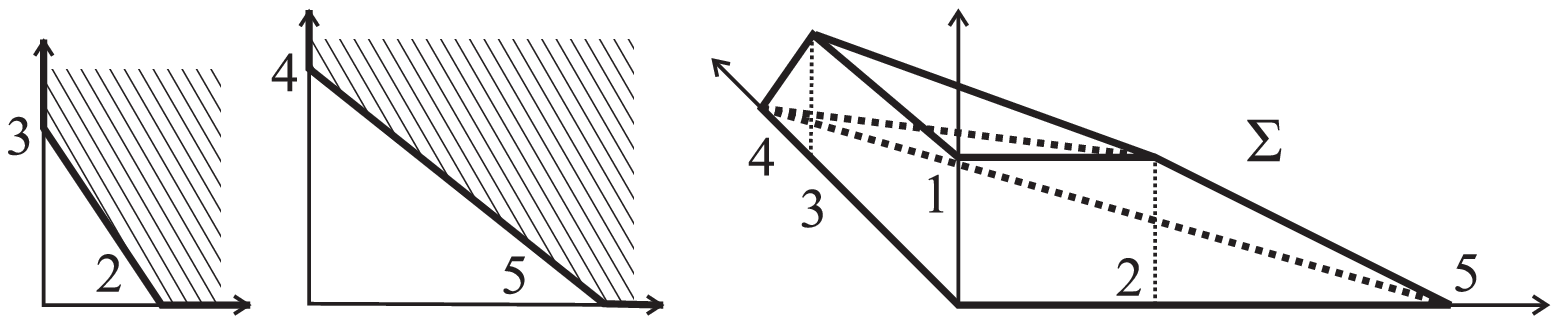}
\end{center}
\end{exa}

Let $N_1,\ldots,N_k,M_1,\ldots,M_m$ be polyhedra in $\R^m_+$, such
that the differences $\R^m_+\setminus N_i$ are bounded. Denote
the mixed volume of pairs $$\bigl(\{0\}\times\R^m_+
,\;\;\{0\}\times N_i\bigr),\;\;i=1,\ldots,k, \mbox{ and } $$ $$
\bigl(M_j*N_1*\ldots*N_k,\;\;M_j*N_1*\ldots*N_k\bigr),\;\;j=1,\ldots,m,$$
in $\R^k\oplus\R^m$ by $\res_{eg}(N_1,\ldots,N_k;M_1,\ldots,M_m)$.
\begin{theor}[\cite{E2}, \cite{E3}]
\label{finoka1} Let $f_1,\ldots,f_k,w_1,\ldots,w_n$ be germs of
holomorphic functions on $(\C^n,0)$, and suppose that their Newton
polyhedra intersect all coordinate axes. If the leading
coefficients of these functions satisfy a certain condition of
general position, then the Gusein-Zade--Ebeling index of the
1-form $w_1 dx_1 + \ldots + w_n dx_n$ on the isolated singularity
of the complete intersection $f_1=\ldots=f_k=0$ makes sense and
equals
$$\sum_{\mathcal{J}=\{i_1,\ldots,i_m\}\subset\{1,\ldots,n\},\atop
\mathcal{J}\ne\emptyset} (-1)^{n-m} (m+k)!
\res_{eg}(\Delta_{f_1}^\mathcal{J},\ldots,\Delta_{f_k}^\mathcal{J};\Delta_{x_{i_1}w_{i_1}}^\mathcal{J},\ldots,\Delta_{x_{i_m}w_{i_m}}^\mathcal{J}).$$
For arbitrary leading coefficients, the index is not smaller than
the expected number, or is not defined.
\end{theor}
\begin{rem} See Theorem \ref{matrzeta}(3) for a generalization to
determinantal singularities under certain additional assumptions.
We do not explicitly formulate the condition of general position
for the leading coefficients, because it is a cumbersome special
case of Definition \ref{defgenpos}. One can reconstruct this
condition, tracing back the reduction of this theorem to Corollary
\ref{volveryconv}.
\end{rem}

The proof is given in Subsection \ref{ss3proof}. In particular,
we can now simplify M.~Oka's formula for the Milnor number of a
complete intersection, formulating it in the language of mixed
volumes of pairs. Consider pairs of polyhedra $A_1,\ldots,A_m$ in
$\R^n$ that have the same support cone $\Gamma$.
\begin{defin}\label{polyseries} We define the value of the power
series \newline
$\sum\limits_{(a_1,\ldots,a_m)\in\Z^n_+}c_{a_1,\ldots,a_m}A_1^{a_1}\cdot\ldots\cdot
A_m^{a_m}$ as the sum $$\sum\limits_{a_1+\ldots+a_m=n}
c_{a_1,\ldots,a_m}\Vol_{\Gamma}(\underbrace{A_1,\ldots,A_1}_{a_1},\ldots,\underbrace{A_m,\ldots,A_m}_{a_m}),$$
and define the \textit{value of a rational function of pairs of
polyhedra} as the value of its power series expansion at the
origin.
\end{defin} For polyhedra $N_0,\ldots,N_k\subset\R^m_+$, such that the
differences $\R^m_+\setminus N_j$ are bounded, we denote the
value of the rational function
$m!\prod\limits_{i=0}^k\frac{(\R^m_+,N_i)}{1+(\R^m_+,N_i)}$ by
$\mu_m(N_0,\ldots,N_k)$.
\begin{sledst}[\cite{oka}] \label{oka1f}
Let $f_0,\ldots,f_k$ be germs of holomorphic functions on
$(\C^n,0)$, such that the differences
$\R^n_+\setminus\Delta_{f_j}$ are bounded. If the leading
coefficients satisfy a certain condition of general position (see
\cite{oka}, or \cite{E1} for a somewhat milder condition), then
the equations $f_0=\ldots=f_k=0$ define an isolated singularity
of a complete intersection, and its Milnor number equals
$$(-1)^{n-k-1} \sum\limits_{\mathcal{J}\subset\{1,\ldots,n\},\;
\mathcal{J}\ne\emptyset}
\mu_{|\mathcal{J}|}(\Delta_{f_0}^\mathcal{J},\ldots,\Delta_{f_k}^\mathcal{J})+(-1)^{n-k}.$$
\end{sledst}
One way to prove this formula is to simplify the answer given in
\cite{oka} by means of Assertion \ref{volsymm}. Another
(independent) way is to compute the Gusein-Zade--Ebeling index of
the 1-form $df_0$ on the complete intersection $f_1=\ldots=f_k=0$
(Theorem \ref{finoka1}), which equals the sum of the Milnor
numbers of the complete intersections $f_0=\ldots=f_k=0$ and
$f_1=\ldots=f_k=0$. The Oka formula then follows by induction on
$k$.

\section{Preliminaries from toric geometry}\label{ss2}

\subsection{Toric varieties} \label{sstor} We recall some basic facts about
smooth toric varieties and introduce corresponding notation (see
details in \cite{danil}, or, more elementary for the smooth case,
in \cite{varch2}). A closed cone in $\R^N$ with the vertex at the
origin is said to be \textit{simple} if it is generated by a part
of a basis of the lattice $\Z^N$. A collection of simple cones is
called a \textit{simple fan} if it is closed with respect to
taking intersections and faces of its cones and satisfies the
following condition: the intersection of each two cones is a face
of both of them. \textit{The support set} $|\Gamma|$ of a fan
$\Gamma$ is defined as the union of its cones.

Simple fans in $\R^N$ are in one-to-one correspondence with
$N$-dimensional smooth toric varieties (i.e. $N$-dimensional
smooth algebraic varieties with an action of a complex torus
$(\C\setminus\{0\})^N$, such that one of the orbits of this
action is everywhere dense -- this orbit is called \textit{the
maximal torus}). We denote the toric variety corresponding to a
fan $\Gamma$ by $\T^{\Gamma}$. The $k$-dimensional orbits of the
variety $\T^{\Gamma}$ are in one-to-one correspondence with the
codimension $k$ cones of the fan $\Gamma$, and adjacent orbits
correspond to adjacent cones. We denote the primitive generator
of a 1-dimensional cone, corresponding to a codimension 1 orbit
$R$, by $\gamma_R$.

Let $\Gamma$ be a simple fan with a convex support set in
$(\R^N)^*$. The polyhedron $\Delta\subset\R^N$ is said to be
\textit{compatible} with $\Gamma$, if its support function
$\Delta(\mathbf{\cdot})$ is linear on every cone of $\Gamma$ and
is not defined outside of $|\Gamma|$. If a polyhedron $\Delta$ is
compatible with $\Gamma$, then there exists a unique ample line
bundle $\mathcal{B}_{\Delta}$ on $\T^{\Gamma}$ equipped with a
meromorphic section $s_{\Delta}$, such that the divisor of zeros
and poles of $s_{\Delta}$ equals $-\sum_{R} \Delta(\gamma_R) R$, 
where $R$ runs over all codimension 1 orbits of $\T^{\Gamma}$.
This correspondence, which assigns the pair
$(\mathcal{B}_{\Delta},s_{\Delta})$ to the polyhedron $\Delta$,
is an isomorphism between the semigroup of convex integer
polyhedra compatible with $\Gamma$, and the semigroup of pairs
$(\mathcal{B},s)$, where $\mathcal{B}$ is a very ample line bundle
on $\T^{\Gamma}$, and $s$ is its meromorphic section with zeroes
and poles outside of the maximal torus.

\textit{The compact part} $\T^{\Int\Gamma}$ of the variety
$\T^{\Gamma}$ is defined as the union of all precompact orbits
(i.e. the orbits, corresponding to the cones in the interior of
$|\Gamma|$). \textit{The dual cone} of the support cone
$|\Gamma|\subset(\R^N)^*$ (i.e. the set of all points in $\R^N$
at which the values of all covectors from $|\Gamma|$ are
non-negative) is denoted by $\Gamma^*\subset\R^N$. Here is a
simple example with $\T^{\Gamma}=\CP^1\times\C^1$ and
$\T^{\Int\Gamma}=\CP^1\times\{0\}$ (we only draw the real part of
the toric variety, of course):
\begin{center} \noindent\includegraphics[width=14cm]{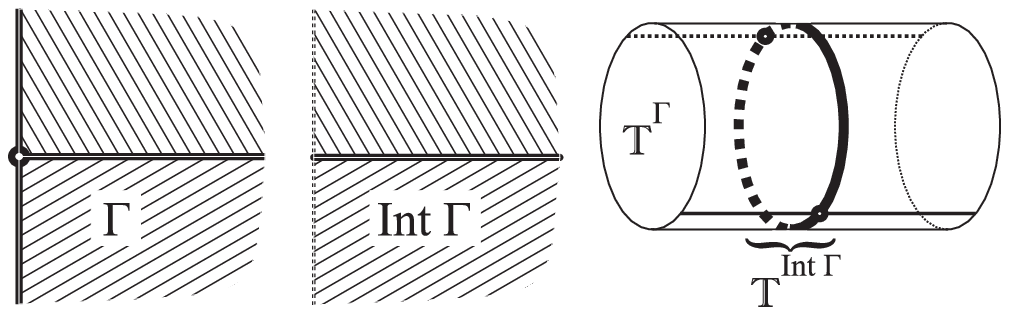}
\end{center}
\begin{defin}[\cite{E3}] \label{defnewtsect}
If a germ $s$ of a meromorphic section of the line bundle
$\mathcal{B}_{\Delta}$ on the pair
$(\T^{\Gamma},\T^{\Int\Gamma})$ has no poles outside of
$\T^{\Int\Gamma}$, then the quotient $s/s_{\Delta}$ defines a
function on $\T^{\Gamma}$ that can be represented as a Laurent
series $\sum_{a\in A} c_a t^a,\; c_a\ne 0,\; A\subset\Z^N$ on the
maximal torus $(\C\setminus\{0\})^N$. In this case:

The \textit{Newton polyhedron} $\Delta_s$ of the germ $s$ is
defined as the convex hull of the set $A+\Gamma^*$;

A coefficient $c_a$ is said to be \textit{leading}, if $a$ is
contained in a bounded face of $\Delta_s$;

For a covector $\gamma$ in the interior of $|\Gamma|$, the lowest
order non-zero $\gamma$-homogeneous component of the series
$\sum_{a\in A} c_a t^a$ is called the
$\gamma$-\textit{truncation} of the section $s$, and is denoted
by $s^{\gamma}$. It is a Laurent polynomial on the maximal torus
$(\C\setminus\{0\})^N$.

The Newton polyhedron of $s$ is said to be not defined, if $s$ has
poles outside of the compact part $\T^{\Int\Gamma}$.
\end{defin}
\begin{rem} In the latter case, we still can represent $s$ as a
quotient of two holomorphic sections $s_1/s_2$, and define the
Newton polyhedron of $s$ as the Minkowski difference
$\Delta_{s_1}-\Delta_{s_2}$, which is a virtual polyhedron. In
this way we can extend most of computations below to arbitrary
meromorphic sections. We do not take this way here, because we do
not need it for our applications in Section \ref{ss1}.
\end{rem}
\begin{exa} \label{exatoric} If $\Delta$ is as shown on the picture
below, then the corresponding line bundle $\mathcal{B}_{\Delta}$
on the toric variety $\T^{\Gamma}=\CP^1\times\C^1$ (see the above
picture) is the pullback of the line bundle $\mathcal{O}(1)$ on
the first factor $\CP^1$. Denoting the standard coordinates on the
factors $\CP^1$ and $\C^1$ by $\lambda:\mu$ and $x$ respectively,
the section $s_{\Delta}$ of the bundle $\mathcal{B}_{\Delta}$ is
the pullback of the section $\mu$ of the bundle $\mathcal{O}(1)$.
The Newton polyhedron of a section $s=\lambda (ax^p+\ldots)+\mu
(bx^q+\ldots)$ of $\mathcal{B}_{\Delta}$ is shown on the picture
below (dots stand for higher order terms), and the leading
coefficients of $s$ are $a$ and $b$.
\end{exa}
\begin{center} \noindent\includegraphics[width=10cm]{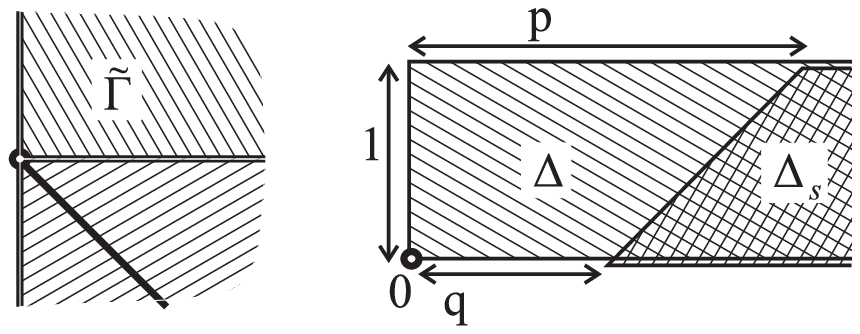}
\end{center}

\subsection{Relative Kouchnirenko-Bernstein formula} The classical
Kouchnirenko-Bernstein formula is about the intersection number
of algebraic hypersurfaces in a complex torus. We need this
formula and the notion of intersection number in a somewhat more
general setting.

\textbf{Intersection numbers.}\label{sskb} We recall the
definition of varieties with multiplicities and their
intersection numbers in the generality that we need (see details
in \cite{fulton}). Let $M$ be a smooth oriented $N$-dimensional
manifold.
\begin{defin}[\cite{E3}]
\textit{A $k$-dimensional cycle on $M$ with the support set $K$}
is defined as a pair $(K,\alpha)$, where $K\subset M$ is a closed
subset, and $\alpha\in
H_k(K\cup\{\infty\},\{\infty\}\,;\;\Q)=H^{N-k}(M,M\setminus
K\,;\;\Q)$.
\end{defin}
Definitions and basic properties of the sum
$(K_1,\alpha_1)+(K_2,\alpha_2)=(K_1\cup K_2, \alpha_1+\alpha_2)$,
the intersection $(K_1,\alpha_1)\cap(K_2,\alpha_2)=(K_1\cap K_2,
\alpha_1\smile\alpha_2)$, the direct image
$f_*(K,\alpha)=(f(K),f_*\alpha)$ under a proper map $f$, and the
inverse image $f^*(K,\alpha)=(f^{(-1)}(K),f^*\alpha)$ under an
arbitrary continuous map $f$ are the same as for homology and
cohomology cycles. If $\alpha$ is the fundamental homology cycle
of a closed irreducible analytical set $K$, then the cycle
$(K,\alpha)$ is denoted by $[K]$ and is called \textit{the
fundamental cycle} of $K$. The divisor of zeroes and poles of a
meromorphic section $w$ of a line bundle $E\to M$ is denoted by
$[w]$. More generally, a continuous section $v$ of a vector
bundle $F\to M$ also corresponds to a certain cycle in $M$, which
is denoted by $[v]$ and is defined as the intersection of the
fundamental cycles of the graphs of the section $v$ and the zero
section in the total space $F$.
\begin{defin}
\textit{The index} $\ind s$ of a 0-dimensional cycle
$s=(K,\alpha)$ with a compact support set $K$ is defined as its
image $\pt_*(\alpha)\in H_0(*\,;\;\Z)=\Z$ under the proper map
$\pt:K\to *$.
\end{defin}
Let $s_1,\ldots,s_k$ be cycles on $M$ such that the sum of their
codimensions equals $N$ and the intersection $K$ of their support
sets is compact (though not necessary 0-dimensional). Then their
\textit{intersection number} $\ind (s_1\cap\ldots\cap s_k)$ makes
sense.

Germs of cycles and indices of their intersections are defined in
the same way. In particular, the index $\ind [A]\cap[B]$ for
germs of analytical sets $A$ and $B$ of complementary dimension
is the intersection number of these germs, the index $\ind [w]$
for the germ of a section $w$ of a rank $N$ vector bundle is the
Poincare-Hopf index of this section.

\textbf{Relative Kouchnirenko-Bernstein formula} Let
$\Delta_1,\ldots,\Delta_I\subset\R^I$ be integer polyhedra
compatible with a simple fan $\Gamma$ in $(\R^I)^*$. Let $s_i,\,
i=1,\ldots,I$ be germs of meromorphic sections of the line bundles
$\mathcal{B}_{\Delta_i}$ on the pair $(\T^{\Gamma},
\T^{\Int\Gamma})$, and suppose that their Newton polyhedra
$\Delta_{s_i}$ are defined (see Definition \ref{defnewtsect}). In
this subsection, we compute the intersection number $\ind
\bigl([s_1]\cap\ldots\cap[s_I]\bigr)$ in terms of polyhedra
$\Delta_i$ and $\Delta_{s_i}$.
\begin{defin}\label{defgenpos0} Leading coefficients of
$s_1,\ldots,s_I$ are said to be \textit{in general position}, if,
for every covector $\gamma$ in the interior of $|\Gamma|$, the
polynomial equations $s_1^{\gamma}=\ldots=s_I^{\gamma}=0$ have no
common roots in $\CC^I$.
\end{defin}
\begin{theor}[Relative Bernstein formula,\cite{E3}] \label{relbernst}
Suppose that, in the above notation, the differences
$\Delta_i\setminus\Delta_{s_i}$ and
$\Delta_{s_i}\setminus\Delta_i, i=1,\ldots,I$ are bounded. If the
leading coefficients of the sections $s_i$ are in general
position, then the intersection number $\ind
([s_1]\cap\ldots\cap[s_I])$ equals $I!$ times the mixed volume of
the pairs $(\Delta_i,\Delta_{s_i})$. If the leading coefficients
are not in general position, then the intersection number is
greater or is not defined.
\end{theor}

\begin{exa}\label{exabernst2} Adopting notation of Example \ref{exatoric} and taking
Example \ref{exavol} into account, the intersection number of the
divisors \begin{center}$\lambda (ax^p+\ldots)+\mu (bx^q+\ldots)=0$
and $\lambda (cx^r+\ldots)+\mu (dx^t+\ldots)=0$\end{center} on the
germ of the variety $\CP^1\times(\C^1,0)$ is greater or equal to
$\min(p+t,q+r)$, and the equality is always attained unless
$p+t=q+r$ and $\det\left(a\; b\atop c\; d\right)=0$. To verify
this answer, we can notice that the desired intersection number
is equal to the multiplicity of zero of the determinant
$$\det\left(ax^p+\ldots\quad bx^q+\ldots\atop cx^r+\ldots\quad
dx^t+\ldots\right) \mbox{ at }  x=0 . $$
\end{exa}

Two proofs of Theorem \ref{relbernst} are given at the end of
this section. The first one (suggested by A.~G.~Khovanskii) is a
very transparent reduction to the global case and gives Assertion
\ref{volexpr2} as a byproduct. The second argument is close to the
original D.~Bernstein's proof and also leads to Assertion
\ref{volsymm} and Corollary \ref{volveryconv}, on which the proof
of Theorem \ref{finoka1} is based. Before proving this theorem,
we explain how to omit the redundant assumption of bounded
differences in Theorem \ref{relbernst}.

\subsection{Stable version of Kouchnirenko-Bernstein formula}\label{ssrelkb} 
As above, let $\Delta_1,\ldots,\Delta_I\subset\R^I$ be integer
polyhedra compatible with a simple fan $\Gamma$ in $(\R^I)^*$.
Let $s_i,\, i=1,\ldots,I$ be germs of meromorphic sections of the
line bundles $\mathcal{B}_{\Delta_i}$ on the pair $(\T^{\Gamma},
\T^{\Int\Gamma})$, and suppose that their Newton polyhedra
$\Delta_{s_i}$ are defined (see Definition \ref{defnewtsect}).

\begin{defin} \label{defgenpos}
Leading coefficients of the sections $s_1,\ldots,s_I$ are said to
be \textit{in strongly general position}, if, for every pair of
covectors $\gamma$ from the boundary of $|\Gamma|$ and $\delta$
from the interior of $|\Gamma|$, the system of polynomial
equations $(s_i^{\gamma})^{\delta}=0,\;
i\in\{i\;|\;\Delta_i(\gamma)=\Delta_{s_i}(\gamma)\}$, has no
solutions in $\CC^I$.
\end{defin}
If the differences of $\Delta_i$ and $\Delta_{s_i}$ are bounded,
then the strongly general position is equivalent to the general
position in the sense of Definition \ref{defgenpos0}. It is
stronger in general, but still occurs for generic leading
coefficients under a certain mild restriction on the Newton
polyhedra (see Corollary \ref{convstab}).

Choose a covector $\gamma_0$ in the interior of $|\Gamma|$, and,
for every positive number $M$, denote the convex hull of the union
$\Delta_{s_i}\cup \bigl(\Delta_i\cap \{\gamma_0\geqslant
M\}\bigr)$ by $H_{i,M}$. If the value of the mixed volume of pairs
$(\Delta_1,H_{1,M}),\ldots,(\Delta_I,H_{I,M})$ does not depend on
$M$ for large $M$, then this value is called the \textit{stable
mixed volume} of pairs
$(\Delta_1,\Delta_{s_1}),\ldots,(\Delta_I,\Delta_{s_I})$.

\begin{exa}
If polygons $A$ and $B$ are as shown below, then the stable mixed
area of the pairs $(\R^2_+,A)$ and $(\R^2_+,B)$ equals
$\frac{ps+rq-pr}{2}$ plus the proper mixed area of pairs
$\bigl(\R^2_++(p,0),A\bigr)$ and $\bigl(\R^2_++(0,r),B\bigr)$.
\end{exa}
\begin{center} \noindent\includegraphics[width=8cm]{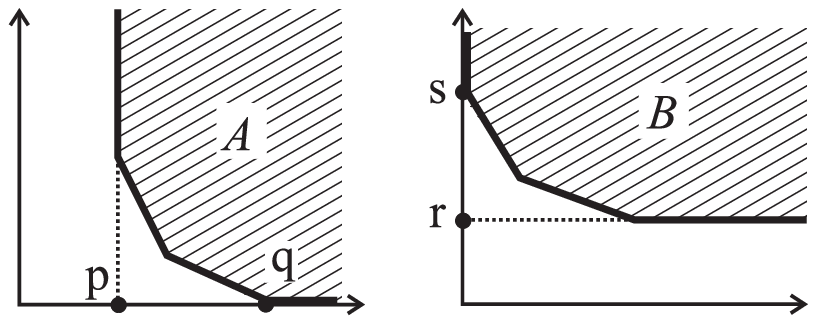}
\end{center}

\begin{theor}\label{relbernstab} If, in the above notation, leading
coefficients of the germs $s_1,\ldots,s_I$ are in strongly general
position, then the intersection index
$\ind[s_1]\cap\ldots\cap[s_I]$ equals $I!$ times the stable mixed
volume of pairs
$(\Delta_1,\Delta_{s_1}),\ldots,(\Delta_I,\Delta_{s_I})$ (in
particular, both of them exist). If leading coefficients of the
germs $s_1,\ldots,s_I$ are not in strongly general position, then
the intersection number $\ind[s_1]\cap\ldots\cap[s_I]$ is greater
or not defined.
\end{theor}
\textsc{Proof.} We can assume without loss of generality that the
polyhedron $\Delta_{i,M}=\Delta_i\cap \{\gamma_0\geqslant M\}$
has integer vertices for every integer $M$. For every
$i=1,\ldots,I$, pick $s_{i,M}$ such that $\Delta_{i,M}$ is its
Newton polyhedron, and consider germs $s_1+t s_{1,M},\ldots,s_I+t
s_{I,M}$ on the pair $\bigl(\T^{\Gamma}\times\C^1,\;
\T^{\Int(\Gamma)}\times\{0\}\bigr)$, where $t$ is the coordinate
on $\C^1$.

If leading coefficients of the collection $s_1,\ldots,s_I$ are in
strongly general position, then

1) one can readily verify that the set $\{s_1+t
s_{1,M}=\ldots=s_I+t s_{I,M}=0\}\subset\T^{\Gamma}\times\C^1$ is
contained in $\T^{\Int(\Gamma)}\times\C^1$ for a large $M$. Thus,
the desired $\ind[s_1]\cap\ldots\cap[s_I]$ equals
$\ind[s_1+\varepsilon s_{1,M}]\cap\ldots\cap[s_I+\varepsilon
s_{I,M}]$ for a small constant $\varepsilon\ne 0$.

2) leading coefficients of the collection $s_1+\varepsilon
s_{1,M},\ldots,s_I+\varepsilon s_{I,M}$ are in general position
for a small constant $\varepsilon\ne 0$. Thus, by Theorem
\ref{relbernst}, $\ind[s_1+\varepsilon
s_{1,M}]\cap\ldots\cap[s_I+\varepsilon s_{1,M}]$ equals $I!$ times
the mixed volume of pairs
$(\Delta_1,\Delta_{s_1})_M,\ldots,(\Delta_I,\Delta_{s_I})_M$.

The first proposition of the theorem follows by (1) and (2).

If leading coefficients of the collection $s_1,\ldots,s_I$ are
not in strongly general position, then we can choose
$s_{1,M},\ldots,s_{I,M}$, such that leading coefficients of the
collection $s_1+\varepsilon s_{1,M},\ldots,s_I+\varepsilon
s_{I,M}$ are not in general position for a small constant
$\varepsilon\ne 0$. Then
$\ind[s_1]\cap\ldots\cap[s_I]\geqslant\ind[s_1+\varepsilon
s_{1,M}]\cap\ldots\cap[s_I+\varepsilon s_{I,M}]$, and the latter
index is greater than $I!$ times the mixed volume of pairs
$(\Delta_1,\Delta_{s_1})_M,\ldots,(\Delta_I,\Delta_{s_I})_M$ by
Theorem \ref{relbernst}, which implies the second proposition of
the theorem. $\Box$

\begin{sledst}\label{convstab} The following three conditions are equivalent:

1) The stable mixed volume of the pairs
$(\Delta_1,\Delta_{s_1}),\ldots,(\Delta_I,\Delta_{s_I})$ exists.

2) Generic leading coefficients are in strongly general position
for sections $s_1,\ldots,s_I$.

3) The collection $\Delta_1,\ldots,\Delta_I,$
$\Delta_{s_1},\ldots,\Delta_{s_I}$ is convenient in the following
sense.
\end{sledst}

\begin{defin} \label{defconv}
In the above notation, the collection of polyhedra
$\Delta_1,\ldots,\Delta_I,$ $\Delta_{s_1},\ldots,\Delta_{s_I}$ is
said to be \textit{convenient}, if, for every pair of covectors
$\gamma$ from the boundary of $|\Gamma|$ and $\delta$ from the
interior of $|\Gamma|$, there exists a set
$I_{\gamma,\delta}\subset \{1,\ldots,I\}$ such that
\newline 1) $\dim \sum_{i\in I_{\gamma,\delta}} (\Delta_{s_i}^{\gamma})^{\delta} \leqslant |I_{\gamma,\delta}|$, \; 2)
$\Delta_{s_i}(\gamma)=\Delta_i(\gamma)$ for every $i\in
I_{\gamma,\delta}$.

\vspace{0.3cm} A convenient collection of polyhedra is said to be
\textit{very convenient} if $\dim \sum_{i\in I_{\gamma,\delta}}
(\Delta_{s_i}^{\gamma})^{\delta} < |I_{\gamma,\delta}|$ for
$\gamma\ne 0$.

\vspace{0.3cm} A very convenient collection of polyhedra is said
to be \textit{cone-convenient} if $\Delta_1$ is a cone and
$I_{\gamma,\delta}\subset \{2,\ldots,n\}$ for every $\gamma\ne 0$.
\end{defin}

\textsc{Proof of Corollary \ref{convstab}.} (2) and (3) are
obviously equivalent. If (2) is satisfied, then (1) is satisfied
by Theorem \ref{relbernstab}. If (3) is not satisfied, then, in
the notation of the proof of Theorem \ref{relbernstab}, we can
choose $s_1,\ldots,s_I$, $s_{1,M},\ldots,s_{I,M}$ and
$s_{1,M+1},\ldots,s_{I,M+1}$ in such a way that leading
coefficients of $s_1+\varepsilon
s_{1,M}+s_{1,M+1},\ldots,s_I+\varepsilon s_{I,M}+s_{I,M+1}$ are
not in general position for a small $\varepsilon\ne 0$, but are in
general position for $\varepsilon=0$. Denote $I!$ times the mixed
volume of the pairs
$(\Delta_1,\Delta_{s_1})_M,\ldots,(\Delta_I,\Delta_{s_I})_M$ by
$V_M$. Then $V_M<\ind[s_1+\varepsilon
s_{1,M}+s_{1,M+1}]\cap\ldots\cap[s_I+\varepsilon
s_{I,M}+s_{I,M+1}]\leqslant\ind[s_1+s_{1,M+1}]\cap\ldots\cap[s_I+s_{I,M+1}]=V_{M+1}$
by Theorem \ref{relbernst}, and (1) is not satisfied. $\Box$

\subsection{Relative Kouchnirenko-Bernstein-Kho\-van\-skii formula}\label{sskbk}
In the assumptions of Theorem \ref{relbernst}, suppose that the
sections $s_1,\ldots,s_k,\; k<I$, are holomorphic (i.e.
$\Delta_{s_i}\subset\Delta_i$ for $i\leqslant k$), and the first
line bundle $\mathcal{I}_{\Delta_1}$ is trivial, i.e. $\Delta_1$
is a cone with the vertex at the origin $0\in\R^I$. The relative
version of the Bernstein-Khovanskii formula computes the Euler
characteristic of the Milnor fiber of the function $s_1$ on the
complete intersection $s_2=\ldots=s_k=0$, in terms of the Newton
polyhedra of the sections $s_1,\ldots,s_k$. At the end of this
subsection, we also explain how to drop the assumption of
triviality of $\mathcal{I}_{\Delta_1}$.

To define the Milnor fiber of $s_1$, it is convenient to fix a
family of neighborhoods for the compact part of the toric variety
$\T^{\Gamma}$. For instance, choose an integer point $a_i$ in
every infinite edge of $\Delta_1$, and define the neighborhood
$B_{\varepsilon}$ of the compact part of $\T^{\Gamma}$ as the
closure of set of all $x\in\CC^I$ such that $\sum_i |x^{a_i}| <
\varepsilon$.
\begin{defin} The Milnor fiber of the function $s_1$ on the
complete intersection $s_2=\ldots=s_k=0$ is the manifold
$\{s_1-\delta=s_2=\ldots=s_k=0\}\cap B_{\varepsilon}$ for
$|\delta|\ll\varepsilon\ll 1$.
\end{defin}
\begin{defin} The leading coefficients of $s_1,\ldots,s_k,\; k<I$, are
said to be in general position, if, for every $\gamma$ in the
interior of $|\Gamma|$, the systems of equations
$s_1^{\gamma}=\ldots=s_k^{\gamma}=0$ and
$s_2^{\gamma}=\ldots=s_k^{\gamma}=0$ define nonsingular varieties
in $\CC^I$.
\end{defin}
\begin{defin} Faces $A_1,\ldots,A_k$ of polyhedra
$\Delta_1,\ldots,\Delta_k$ in $\R^I$ are said to be compatible, if
$A_i=\Delta_i^{\gamma},\; i=1,\ldots,k$, for some linear function
$\gamma$ on $\R^I$.
\end{defin}

For compatible unbounded faces
$A_1\subset\Delta_1,\ldots,A_k\subset\Delta_k$, denote the
dimension of $A_1+\ldots+A_k$ by $q$. Then, up to a shift, the
pairs $(A_i,A_i\cap\Delta_{s_i})$ are contained in the same
rational $q$-dimensional subspace of $\R^N$; we denote the
corresponding shifted pairs by $\widetilde A_1,\ldots,\widetilde
A_k$. In the notation, introduced prior to Corollary \ref{oka1f},
denote the number $q!\frac{\widetilde
A_1\cdot\ldots\cdot\widetilde A_k}{(1+\widetilde
A_1)\cdot\ldots\cdot(1+\widetilde A_k)}$ by
$\chi_{A_1,\ldots,A_k}$.

\begin{theor} \label{relbkh} In the above assumptions, the Euler characteristic of the Milnor fiber of $s_1$
on the complete intersection $\{s_2=\ldots=s_n=0\}$ equals the
sum of $\chi_{A_1,\ldots,A_k}$ over all collections
$A_1\subset\Delta_1,\ldots,A_k\subset\Delta_k$ of compatible
unbounded faces, provided that the leading coefficients of
$s_1,\ldots,s_k$ are in general position.
\end{theor}

This is proved in \cite{oka} for a regular affine toric variety
(based on the idea of \cite{varch}), and in \cite{takeuchi} for
an arbitrary affine toric variety (in a more up to date
language). Both arguments can be easily applied to an arbitrary
(not necessary affine) toric variety, and also provide a formula
for the $\zeta$-function of monodromy of the function $s_1$.

However, since we restrict our consideration by the Milnor number
in this paper, we prefer to give a much simpler proof by
reduction to the global Bernstein-Khovanskii formula. It is based
on the same Khovanskii's idea, as the first proof of Theorem
\ref{relbernst}.

\textsc{Proof (\cite{E5}).} If the leading coefficients are in
general position, then topology of the Milnor fiber only depends
on the Newton polyhedra of $s_1,\ldots,s_k$, and we can assume
without loss of generality that $s_i=s_{\Delta_i}\cdot\tilde
s_i$, where $\tilde s_1,\ldots,\tilde s_k$ are Laurent
polynomials on $\CC^I$ and satisfy the condition of general
position of \cite{pkhovvol}.

Pick an orbit $T$ of $\T^{\Gamma}$, and consider the
corresponding faces $\widetilde D_i$ and $D_i$ of the polyhedron
$\Delta_i$ and the Newton polyhedron of $\tilde s_i$ (i.e. the
faces $\Delta^{\gamma}_i$ and $\Delta^{\gamma}_{\tilde s_i}$ for
a covector $\gamma$ in the relative interior of the cone,
corresponding to the orbit $T$ in the fan $\Gamma$). Up to a
shift, all $\widetilde D_i$ and $D_i,\; i=1,\ldots,k$, are
contained in the same $(\dim T)$-dimensional subspace of $\R^I$,
hence, their $(\dim T)$-dimensional mixed volumes make sense.

By the global Bernstein-Khovanskii formula, the Euler
characteristics of $T\cap\{\tilde s_1=\ldots=\tilde s_k=0\}$ and
$T\cap\{\tilde s_1-\delta=\tilde s_2=\ldots=\tilde s_k=0\}$ equal
$(\dim T)!\frac{D_1\cdot\ldots\cdot
D_k}{(1+D_1)\cdot\ldots\cdot(1+D_k)}$ and $(\dim
T)!\frac{\widetilde D_1\cdot D_2\cdot\ldots\cdot
D_k}{(1+\widetilde D_1)\cdot(1+D_2)\cdot\ldots\cdot(1+D_k)}$
respectively.

Since the boundary of $B_{\varepsilon}$ subdivides the set
$T\cap\{\tilde s_1-\delta=\tilde s_2=\ldots=\tilde s_k=0\}$ into
two parts, homeomorphic to the set $T\cap\{\tilde
s_1=\ldots=\tilde s_k=0\}$ and the Milnor fiber of $s_1$ on
$T\cap\{s_2=\ldots=s_n=0\}$, the Euler characteristic of the
latter equals $$(\dim T)!\frac{\widetilde D_1\cdot
D_2\cdot\ldots\cdot D_k}{(1+\widetilde
D_1)\cdot(1+D_2)\cdot\ldots\cdot(1+D_k)}-(\dim
T)!\frac{D_1\cdot\ldots\cdot
D_k}{(1+D_1)\cdot\ldots\cdot(1+D_k)}$$ by additivity of Euler
characteristic.

Denote the minimal face of $\Delta_i$, containing $D_i$, by $A_i$.
Then the difference above equals $\chi_{A_1,\ldots,A_k}$ by
Proposition \ref{volexpr2}. Thus, the Euler characteristic of the
Milnor fiber of $s_1$ on $T\cap\{s_2=\ldots=s_n=0\}$ equals
$\chi_{A_1,\ldots,A_k}$.

Summing up these equalities over all orbits
$T\subset\T^{\Gamma}$, corresponding to the cones from the
boundary of $\Gamma$, the statement of the theorem follows by
additivity of Euler characteristic. $\Box$

Finally, we formulate a more general version of this theorem,
with no assumptions on the triviality of the first line bundle.
We omit the proof, since it is essentially the same, and we do
not need this statement for our applications to determinantal
singularities.

In the assumptions of Theorem \ref{relbernst}, suppose that the
sections $s_1,\ldots,s_k$, $k<I$, are holomorphic (i.e.
$\Delta_{s_i}\subset\Delta_i$ for $i\leqslant k$), and pick a
holomorphic section $t_i$ of the line bundle
$\mathcal{I}_{\Delta_i}$ in the neighborhood $B_{\varepsilon}$
for a small $\varepsilon$. Varieties
$B_{\varepsilon}\cap\{s_1-t_1=\ldots=s_k-t_k=0\}\setminus\T^{\Int\Gamma}$
are diffeomorphic to each other for almost all collections of
small sections $(t_1,\ldots,t_k)$. Such variety is called the
Milnor fiber of the complete intersection $\{s_1=\ldots=s_n=0\}$.

\begin{theor} In the above assumptions, the Euler characteristic of the Milnor fiber of
the complete intersection $\{s_1=\ldots=s_n=0\}$ equals the sum
of $\chi_{A_1,\ldots,A_k}$ over all collections
$A_1\subset\Delta_1,\ldots,A_k\subset\Delta_k$ of compatible
unbounded faces, provided that the leading coefficients of
$s_1,\ldots,s_k$ are in general position.
\end{theor}

\subsection{Proof of Theorem \ref{relbernst}}\label{ss2proof}

\textbf{Preliminary remarks.} The following lemma clarifies the
statement of Theorem \ref{relbernst} and reduces it to the case
of germs $s_1,\ldots,s_I$ with generic leading coefficients.

\begin{lemma}\label{l51}
Under the assumptions of Theorem \ref{relbernst}
\newline 1) the intersection number $\ind
([s_1]\cap\ldots\cap[s_I])$ does not depend on the choice of the
simple fan $\Gamma$;
\newline 2) intersection number $\ind
([s_1]\cap\ldots\cap[s_I])$ makes sense if leading coefficients of
$s_1,\ldots,s_I$ are in general position.
\newline 3)
Let $\tilde s_i$ be a section of the line bundle
$\mathcal{B}_{\Delta_i}$ such that $\Delta_{\tilde
s_i}=\Delta_{s_i}$, and suppose that leading coefficients of the
sections $s_1,\ldots,s_I$ are in general position. Then $\ind
([\tilde s_1]\cap\ldots\cap[\tilde s_I])\; \geqslant\; \ind
([s_1]\cap\ldots\cap[s_I])$, and the equality takes place if and
only if leading coefficients of the sections $\tilde
s_1,\ldots,\tilde s_I$ are in general position.
\end{lemma}
Note that Theorem \ref{relbernstab} provides a generalization of
these propositions to the case of a convenient collection
$\Delta_1,\ldots,\Delta_I,$ $\Delta_{s_1},\ldots,\Delta_{s_I}$,
but neither of them are valid if the collection of the Newton
polyhedra is not convenient.
\newline
\textsc{Proof.} 1) Suppose that simple fans $\Gamma_1$ and
$\Gamma_2$ are compatible with the polyhedra
$\Delta_1,\ldots,\Delta_I$. If one of them is a subdivision of
the other one, then there exists the natural mapping
$p:\T^{\Gamma_1}\to\T^{\Gamma_2}$ of topological degree 1, such
that $\Delta_{p^*s_i}=\Delta_{s_i}$. Thus, intersection numbers
$\ind [s_1]\cap\ldots\cap[s_I]$ and $\ind
[p^*s_1]\cap\ldots\cap[p^*s_I]$ on the varieties $\T^{\Gamma_1}$
and $\T^{\Gamma_2}$ are equal, and Part 1 follows. In general,
neither of $\Gamma_1$ and $\Gamma_2$ is a subdivision of the
other one, but they always admit a common simple subdivision.
\newline 2) Part 1 allows us to assume that the fan $\Gamma$ is compatible with all
polyhedra $\Delta_i$ and $\Delta_{s_i},\, i=1,\ldots,I$.
Represent every cycle $[s_i]$ as a sum $\alpha_i+\beta_i$, where
the support set $|\alpha_i|$ is contained in the compact part
$\T^{\Int(\Gamma)}$, and the support set $|\beta_i|$ intersects
$\T^{\Int(\Gamma)}$ properly. Then the condition of general
position for leading coefficients of $s_1,\ldots,s_I$ can be
reformulated as follows:
$$|\beta_1|\cap\ldots\cap|\beta_I|=\varnothing.$$
Hence, $[s_1]\cap\ldots\cap [s_I]$ equals the sum of intersections
$\cap_{i\in\mathcal{I}}\alpha_i\cap_{i\notin\mathcal{I}}\beta_i$
over all non-empty sets $\mathcal{I}\subset\{1,\ldots,I\}$. Since
the support set of each of these intersections is compact, their
indices make sense.
\newline 3) Represent every $s_i$ as $\alpha_i+\beta_i$ and $\tilde s_i$ as
$\tilde\alpha_i+\tilde\beta_i$ as above, and note that
$\alpha_i=\tilde\alpha_i$ because $\Delta_{s_i}=\Delta_{\tilde
s_i}$. Thus, taking into account that
$|\beta_1|\cap\ldots\cap|\beta_I|=\varnothing$, we have
$$\ind([\tilde s_1]\cap\ldots\cap [\tilde s_I])=\ind([s_1]\cap\ldots\cap [s_I])+\ind(\tilde\beta_1\cap\ldots\cap\tilde\beta_I).$$
If leading coefficients of the sections $\tilde s_1,\ldots,\tilde
s_I$ are in general position, then
$|\tilde\beta_1|\cap\ldots\cap|\tilde\beta_I|=\varnothing$ and
$\ind[\tilde s_1]\cap\ldots\cap [\tilde s_I]=\ind
[s_1]\cap\ldots\cap [s_I]$, otherwise
$|\tilde\beta_1|\cap\ldots\cap|\tilde\beta_I|\ne\varnothing$. In
the latter case, if
$|\tilde\beta_1|\cap\ldots\cap|\tilde\beta_I|$ is not compact,
then $\ind ([\tilde s_1]\cap\ldots\cap[\tilde s_I])$ is not
defined. Otherwise we have
$\ind\tilde\beta_1\cap\ldots\cap\tilde\beta_I>0$, because, upon a
small perturbation of $\tilde s_1,\ldots,\tilde s_I$, the set
$|\tilde\beta_1|\cap\ldots\cap|\tilde\beta_I|$ is finite and
non-empty. $\Box$

\textbf{Proof of Theorem \ref{relbernst} and Assertion
\ref{volexpr2}.} By Lemma \ref{l51}, Part 3, it is enough to
prove Theorem \ref{relbernst} under the assumption that leading
coefficients of the sections $s_i,\, i=1,\ldots,I$ are in general
position, and the quotient $s_i/s_{\Delta_i}$ defines a Laurent
polynomial on the maximal torus (recall that $s_{\Delta_i}$ is
the distinguished section of the line bundle
$\mathcal{B}_{\Delta_i}$ with no zeros and poles on the maximal
torus, see Subsection \ref{sstor}). Since, under this assumption,
the intersection number $\ind ([s_1]\cap\ldots\cap[s_I])$ is a
symmetric multilinear function of the pairs of polyhedra
$(\Delta_i,\Delta_{s_i})$, it is enough to prove the equality
$\ind
([s_1]\cap\ldots\cap[s_I])=I!\Vol(\Delta_1\setminus\Delta_{s_1})-I!\Vol(\Delta_{s_i}\setminus\Delta_i)$
under the additional assumption that $\Delta_1=\ldots=\Delta_I$,
$\Delta_{s_1}=\ldots=\Delta_{s_I}$. The proof is based on the
following fact. 
\begin{lemma}
Suppose that, under the assumptions of Theorem \ref{relbernst},
$s_i/s_{\Delta_i}=\tilde s_i$, where $\tilde s_1,\ldots,\tilde
s_I$ are Laurent polynomials with (bounded) Newton polyhedra
$\Delta_{\tilde s_i}$, and leading coefficients of $\tilde
s_1,\ldots,\tilde s_I$ are in general position. Then
$$\ind
([s_1]\cap\ldots\cap[s_I])=I!\Vol(\tilde\Delta_1,\ldots,\tilde\Delta_I)-I!\Vol(\Delta_{\tilde
s_1},\ldots,\Delta_{\tilde s_I}),$$ where
$\tilde\Delta_i=\Delta_i\setminus(\Delta_{s_i}\setminus\Delta_{\tilde
s_i})$.
\end{lemma}
\textsc{Proof.} We extend $\Gamma$ to a complete fan
$\tilde\Gamma$, compatible with the polyhedra $\tilde\Delta_i$ and
$\Delta_{\tilde s_i}, i=1,\ldots,I$. The toric variety
$\T^{\tilde\Gamma}$ is a compactification of $\T^{\Gamma}$. The
line bundles $\mathcal{B}_{\Delta_i}$ and their sections $s_i$
extend to the line bundles $\mathcal{B}_{\tilde\Delta_i}$ and
their sections $q_i=s_{\tilde\Delta_i}\cdot\tilde s_i$ on this
compactification. The set of common zeroes of the sections
$q_1,\ldots,q_I$ consists of a compact set $Z\subset\T^{\Int
\Gamma}$ and finitely many points in the maximal torus, which
splits the intersection $[q_1]\cap\ldots\cap[q_I]$ into the sum
$q_Z+q_T$, where $|q_Z|=Z$ and $|q_T|$ is contained in the
maximal torus: $$[q_1]\cap\ldots\cap[q_I]=q_Z+q_T.$$ By the
classical D.~Bernstein's formula, formulated in terms of
intersection numbers of divisors on toric varieties (see
\cite{pkhovvol}), we have \begin{center} $\ind
([q_1]\cap\ldots\cap[q_I])=I!\Vol(\tilde\Delta_1,\ldots,\tilde\Delta_I)$
and $\ind q_T = I!\Vol(\Delta_{\tilde s_1},\ldots,\Delta_{\tilde
s_I})$.\end{center} On the other hand, $\ind q_Z$ equals the
desired intersection number $\ind ([s_1]\cap\ldots\cap[s_I])$ by
the construction of $q_Z$. $\Box$

If we have $\Delta_1=\ldots=\Delta_I$,
$\Delta_{s_1}=\ldots=\Delta_{s_I}$, leading coefficients of the
sections $s_i,\, i=1,\ldots,I$ are in general position, and
$s_i/s_{\Delta_i},\, i=1,\ldots,I$ are Laurent polynomials, then
this lemma implies that $\ind
([s_1]\cap\ldots\cap[s_I])=I!\Vol(\Delta_1\setminus\Delta_{s_1})-I!\Vol(\Delta_{s_i}\setminus\Delta_i)$,
and Theorem \ref{relbernst} follows. Together with Theorem
\ref{relbernst}, this lemma proves Assertion \ref{volexpr2}.
$\square$

\textbf{Proof of Theorem \ref{relbernst} and Assertion
\ref{volsymm} (\cite{E3}).} It is enough to prove Theorem
\ref{relbernst} under the assumption that
$\Delta_{s_i}\subset\Delta_i$ (which means that $s_1,\ldots,s_I$
are holomorphic sections). Indeed, we can represent the
meromorphic section $s_i$ as the quotient of the sections
$s_iq_i$ and $q_i$ of the line bundles
$\mathcal{B}_{\Delta_i+\tilde\Delta_i}$ and
$\mathcal{B}_{\tilde\Delta_i}$ respectively, such that
$\Delta_{s_iq_i}=\Delta_{s_i}+\Delta_{q_i}\subset\Delta_i+\tilde\Delta_i$
and $\Delta_{q_i}\subset\tilde\Delta_i$ (which means that
$s_iq_i$ and $q_i$ are holomorphic). The intersection number of
the divisors of sections is additive with respect to the
operation of tensor multiplication of sections, the mixed volume
of pairs is additive with respect to the operation of Minkowski
summation of pairs, hence the statement of the theorem for the
collection of meromorphic sections $s_1,\ldots,s_I$ follows from
the same statement for all collections of holomorphic sections of
the form $q_1(s_1)^{\alpha_1},\ldots,q_I(s_I)^{\alpha_I}$, where
$(\alpha_1,\ldots,\alpha_I)$ ranges over $\{0,1\}^I$.

First, we compute the desired intersection number in the
following special case (see Subsection \ref{ssmixvol} for notation
related to support functions and faces of convex polyhedra).
\begin{lemma} \label{indvoltriv}
Let $\Delta_1,\ldots,\Delta_I,$
$\Delta_{s_1},\ldots,\Delta_{s_I}$ be a cone-convenient collection
of polyhedra (see Definition \ref{defconv}). If leading
coefficients of the sections $s_i$ are in strongly general
position, then the intersection number $\ind
([s_1]\cap\ldots\cap[s_I])$ equals
$$
(I-1)! \sum\limits_{\gamma\in S\cap\Int|\Gamma|}
(\Delta_{s_1}(\gamma) - \Delta_1(\gamma))\cdot \Vol
(\Delta_{s_2}^{\gamma},\ldots,\Delta_{s_I}^{\gamma}).
$$
(recall that $\Vol$ is the classical mixed volume of bounded
polyhedra, $S$ stands for the set of all primitive covectors in
$(\Z^I)^*$ and $\Int|\Gamma|$ stands for the interior of the cone
$|\Gamma|$; one can readily verify that the right hand side
contains finitely many non-zero terms).
\end{lemma}
\textsc{Proof.} We choose a simple fan $\Gamma$ compatible with
the polyhedra $\Delta_{s_i},\Delta_i,\, i=1,\ldots,I$. The
divisor $[s_i]$ can be represented as
$\alpha_i+\beta_i+\varphi_i$, where $\alpha_i$ is a linear
combination of the closures of precompact codimension 1 orbits in
$\T^{\Gamma}$, $\beta_i$ is a linear combination of the closures
of non-precompact codimension 1 orbits, and the hypersurfaces
$|\varphi_i|, i=1,\ldots,I$ intersect codimension 1 orbits of
$\T^{\Gamma}$ properly.

We note that the desired intersection number $\ind
(\alpha_1+\beta_1+\varphi_1)\cap\ldots\cap(\alpha_I+\beta_i+\varphi_I)$
equals $$\ind\alpha_1\cap\varphi_2\cap\ldots\cap\varphi_I.$$
Indeed, since the line bundle $\mathcal{B}_{\Delta_1}$ is
trivial, the support set of the cycle $[s_1-\varepsilon]\cap
\alpha_j$ is empty, and, in particular, $\ind
[s_1]\cap\ldots\cap\ldots\cap[s_I]=\ind
[s_1-\varepsilon]\cap[s_2]\cap\ldots\cap[s_I]=\ind
[s_1-\varepsilon]\cap(\beta_2+\varphi_2)\cap\ldots\cap(\beta_I+\varphi_I)=\ind
[s_1]\cap(\beta_2+\varphi_2)\cap\ldots\cap(\beta_I+\varphi_I)$.
Since the polyhedra under consideration are cone-convenient and
leading coefficients of $s_1,\ldots,s_I$ are in strongly general
position, the support sets of the intersections
$(\beta_1+\varphi_1)\cap\ldots\cap (\beta_I+\varphi_I)$ and
$\alpha_1\cap\bigcap_{j\in
\mathcal{J}}\varphi_j\cap\bigcap_{j\notin \mathcal{J}}\beta_j,\;
\mathcal{J}\subsetneq\{1,\ldots,I\}$ are empty, and hence $\ind
(\alpha_1+\beta_1+\varphi_1)\cap(\beta_2+\varphi_2)\cap\ldots\cap(\beta_I+\varphi_I)=\ind\alpha_1\cap\varphi_2\cap\ldots\cap\varphi_I$.

We can compute the intersection number
$$\ind\alpha_1\cap\varphi_2\cap\ldots\cap\varphi_I,$$ using the following two equalities.
By D.~Bernstein's theorem, if $R$ is a compact closure of a
codimension 1 orbit of $\T^{\Gamma}$, and $\gamma_R$ is the
primitive generator of the corresponding 1-dimensional cone of
$\Gamma$, then the intersection number $\ind
[R]\cap\varphi_2\cap\ldots\cap\varphi_I$ equals $(I-1)!\Vol
(\Delta_{s_2}^{\gamma_R},\ldots,\Delta_{s_I}^{\gamma_R})$. On the
other hand, $\alpha_1=\sum_{R}
\Bigl(\Delta_{s_1}(\gamma_R)-\Delta_1(\gamma_R)\Bigr)\cdot[R]$,
where $R$ runs over all compact closures of codimension 1 orbits
of $\T^{\Gamma}$.

The second statement of the lemma follows from Lemma \ref{l51},
Part 3. $\Box$

We can reduce the general problem of computation of the
intersection number $\ind ([s_1]\cap\ldots\cap[s_I])$ to the
special case of Lemma \ref{indvoltriv} as follows (perturbation
of divisors that we use in this argument is similar to the
original idea of D.~Bernstein, see \cite{bernst}).

\begin{lemma} \label{indvolgen}
Suppose that the collection of polyhedra
$\Delta_{s_i}\subset\Delta_i\subset\R^I,\; i=1,\ldots,I$ is very
convenient. For arbitrary numbers $a_i\in\N, i=1,\ldots,I$, denote
by $\Sigma_i$ the convex hull of
$((\R^1_++a_i)\times\Delta_i)\cup(\{0\}\times\Delta_{s_i})\subset\R^1\oplus\R^I$,
and denote the point $(1,0,\ldots,0)\in\R^1\oplus\R^I$ by $E$. If
leading coefficients of the sections $s_i$ are in strongly general
position, then the intersection number $\ind
([s_1]\cap\ldots\cap[s_I])$ equals $ I! \sum\limits_{\gamma\in
S\cap\Int(\R^1_+\times|\Gamma|)} (\gamma, E) \Vol
(\Sigma_1^{\gamma},\ldots,\Sigma_I^{\gamma})$ (recall that $\Vol$
is the classical mixed volume of bounded polyhedra, $S$ stands
for the set of all primitive covectors in $(\Z^1\oplus\Z^I)^*$
and $\Int$ stands for the interior; one can readily verify that
the right hand side contains only finitely many non-zero terms).
\end{lemma}
\textsc{Proof.} We pick sections $\tilde s_i$ of the line bundles
$\mathcal{B}_{\Delta_i}$ with the Newton polyhedra
$\Delta_{\tilde s_i}=\Delta_i$ and generic leading coefficients.
Denote by $\varepsilon$ the standard coordinate on the first
factor of the product
$\C^1\times\T^{\Gamma}=\T^{\R^1_+\times\Gamma}$. The desired
intersection number $\ind ([s_1]\cap\ldots\cap[s_I])$ equals the
intersection number
$\ind[\varepsilon]\cap[s_1+\varepsilon^{a_1}\tilde
s_1]\cap\ldots\cap[s_I+\varepsilon^{a_I}\tilde s_I]$ on the germ
of the pair $(\C^1\times\T^{\Gamma},\{0\}\times\T^{\Int\Gamma})$.
This intersection number can be computed by Lemma
\ref{indvoltriv}, because $\varepsilon$ can be considered as a
section of the trivial line bundle on $\C^1\times\T^{\Gamma}$.
$\Box$

In the notation of Lemma \ref{indvolgen}, a bounded face of the
polyhedron $\Sigma_i\subset\R^1\oplus\R^I$ is said to be
\textit{long}, if its projection to the first summand of the
space $\R^1\oplus\R^I$ is not a point. Denote by $L_i$ the set of
all covectors $\gamma$ such that the support face of $\Sigma_i$
with respect to $\gamma$ is long. If the sequence $a_1\ll
a_2\ll\ldots\ll a_I$ is fast-increasing enough, then the sets
$L_i$ do not intersect each other. The answer that Lemma
\ref{indvolgen} gives for $a_1\ll a_2\ll\ldots\ll a_I$ can be
simplified by Lemma \ref{comput1}, and turns out to be equal to
$\sum_{k=1}^I \sum\limits_{\gamma\in S\cap\Int|\Gamma|}
(\Delta_{s_k}(\gamma)-\Delta_k(\gamma))
\Vol(\Delta_1^{\gamma},\ldots,
\Delta_{k-1}^{\gamma},\Delta_{s_{k+1}}^{\gamma},\ldots,\Delta_{s_I}^{\gamma}).$

This expression is a symmetric function of pairs
$(\Delta_i,\Delta_{s_i})$, because the intersection number $\ind
([s_1]\cap\ldots\cap[s_I])$ is a symmetric function of
$s_1,\ldots,s_I$. By Theorem \ref{voldef}, Part 2, the
symmetrization of this expression equals the mixed volume of the
pairs $(\Delta_i,\Delta_{s_i})$, hence this expression itself
equals the mixed volume, which proves Theorem \ref{relbernst} and
also implies Assertion \ref{volsymm}. $\square$

To prove Theorem \ref{finoka1} we need the following version of
Theorem \ref{relbernst}:
\begin{sledst} \label{volveryconv}
Suppose that integer polyhedra
$\Delta_1,\ldots,\Delta_I\subset\R^I$ are compatible with a
simple fan $\Gamma$ in $(\R^I)^*$. Let $s_i, i=1,\ldots,I,$ be
germs of holomorphic sections of the line bundles
$\mathcal{B}_{\Delta_i}$ on the pair $(\T^{\Gamma},
\T^{\Int\Gamma})$, suppose that the collection of polyhedra
$\Delta_1,\ldots,\Delta_I,$ $\Delta_{s_1},\ldots,\Delta_{s_I}$ is
cone-convenient, and the difference
$\Delta_1\setminus\Delta_{s_1}$ is bounded. If leading
coefficients of the sections $s_i$ are in strongly general
position, then the intersection number $\ind
([s_1]\cap\ldots\cap[s_I])$ equals
$$
I!\Vol_{|\Gamma|}\Bigl((\Delta_1,\Delta_{s_1}),(\Delta_{s_2},\Delta_{s_2}),\ldots,(\Delta_{s_I},\Delta_{s_I})\Bigr).
$$
\end{sledst}
\textsc{Proof.} The intersection number can be computed by Lemma
\ref{indvoltriv}. We prove that the answer equals
$I!\Vol_{|\Gamma|}\Bigl((\Delta_1,\Delta_{s_1}),(\Delta_{s_2},\Delta_{s_2}),\ldots,(\Delta_{s_I},\Delta_{s_I})\Bigr)$,
rewriting the latter mixed volume of pairs by Assertion
\ref{volsymm}. $\Box$

\section{Resultantal varieties}\label{ss3}

We introduce resultantal singularities and study their invariants
in terms of Newton polyhedra, which relies upon toric geometry of
Section \ref{ss2} and includes results of Section \ref{ss1} as a
special case.

\subsection{Resultantal varieties}\label{ssresvar}

We denote the monomial $t_1^{a_1}\ldots t_N^{a_N}$ by $t^a$. For
a subset $\Sigma\subset\R^N$, we denote the set of all Laurent
polynomials of the form $\sum\limits_{a\in\Sigma\cap\Z^N} c_a
t^a,\, c_a\in\C$, by $\C[\Sigma]$. We regard them as functions on
the complex torus $(\C\setminus\{0\})^N$.
\begin{defin}[\cite{E2}, \cite{E3}]
\label{torrezcikl} For finite sets $\Sigma_i\subset \Z^N,
i=1,\ldots,I$, the \textit{resultantal variety}
$R(\Sigma_1,\ldots,\Sigma_I)$ is defined as the closure of the set
$$\bigl\{ (g_1,\ldots,g_I)\; |\; g_i\in\C[\Sigma_i],\; \exists
t\in(\C\setminus\{0\})^N \; :$$ $$
g_1(t)=\ldots=g_I(t)=0\bigr\}\subset
\C[\Sigma_1]\oplus\ldots\oplus\C[\Sigma_I].$$
\end{defin}
\begin{defin}[\cite{E2}, \cite{E3}] \label{defressing} Consider a germ of a holomorphic mapping
$f:(\C^n,0)\to(\C[\Sigma_1]\oplus\ldots\oplus\C[\Sigma_I],0)$. The
preimage $f^{(-1)}\bigl(R(\Sigma_1,\ldots,\Sigma_I)\bigr)$ is
called a \textit{resultantal singularity}, if its codimension
equals the codimension of $R(\Sigma_1,\ldots,\Sigma_I)$.
\end{defin}
For instance, a determinantal singularity is a special case of a
resultantal singularity by the following lemma.
\begin{lemma} \label{detviarez} Identify an $I\times k$ matrix
$(w_{i,\, l})$ with the collection of linear functions
$$w_{1,\, l}t_1+\ldots+w_{I-1,\, l}t_{I-1}+w_{I,\, l},\; l=1,\ldots,k.$$ Then the
space of all collections of matrices $\C^{I_1\times
k_1}\oplus\ldots\oplus\C^{I_J\times k_J}$ is identified with the
space
$$\underbrace{\C[\Sigma_1]\oplus\ldots\oplus\C[\Sigma_1]}_{k_1}\oplus\ldots\oplus
\underbrace{\C[\Sigma_J]\oplus\ldots\oplus\C[\Sigma_J]}_{k_J},$$
where $\Sigma_j$ is the set of vertices of the standard
$(I_j-1)$-dimensional simplex in the $j$-th summand of the direct
sum ${\mathbf{Z}}^{I_1-1}\oplus\ldots\oplus{\mathbf{Z}}^{I_J-1}$,
and the set of all collections of degenerate matrices is
identified with the resultantal variety
$$R(\underbrace{\Sigma_1,\ldots,\Sigma_1}_{k_1},\ldots,
\underbrace{\Sigma_J,\ldots,\Sigma_J}_{k_J}).$$
\end{lemma}
\begin{rem} It would be useful to extend well-known relations
between determinantal varieties and determinants to resultantal
varieties and resultants: to prove, for instance, that the ideal
of the resultantal variety $R(\Sigma_1,\ldots,\Sigma_I)$ is
generated by $(\Sigma_{i_1},\ldots,\Sigma_{i_p})$-resultants,
where $(\Sigma_{i_1},\ldots,\Sigma_{i_p})$ runs over all essential
subcollections of codimension 1.
\end{rem}
When studying resultantal singularities, we can impose some
helpful restrictions on the sets $\Sigma_1,\ldots,\Sigma_I$
without loss in generality. We denote the dimension of the convex
hull of $\Sigma_1+\ldots+\Sigma_I$ by
$\dim(\Sigma_1,\ldots,\Sigma_I)$, and the difference
$I-\dim(\Sigma_1,\ldots,\Sigma_I)$ by
$\codim(\Sigma_1,\ldots,\Sigma_I)$.
\begin{lemma}[\cite{sturmf}] \label{sturm1} We have
$$\codim
R(\Sigma_1,\ldots,\Sigma_I)=\max\limits_{\{j_1,\ldots,j_p\}\subset\{1,\ldots,I\}}\codim(\Sigma_{j_1},\ldots,\Sigma_{j_p}).$$
\end{lemma}
See \cite{pkhovvol} and \cite{sturmf} for the special cases of
$\codim(\Sigma_1,\ldots,\Sigma_I)=0$ and $1$ respectively. The
generalization for arbitrary $\codim(\Sigma_1,\ldots,\Sigma_I)$
is quite straightforward, see e.g. \cite{E4}, Theorem 2.12, for
details.
\begin{defin}[\cite{sturmf}]
A collection of sets $\Sigma_i\subset \Z^N, i=1,\ldots,I$, is said
to be \textit{essential}, if
$\codim(\Sigma_{i_1},\ldots,\Sigma_{i_J})<\codim(\Sigma_1,\ldots,\Sigma_I)$
for every subset $\{i_1,\ldots,i_J\}\subsetneq\{1,\ldots,I\}$.
\end{defin}
The following lemma implies that, without loss of generality, we
can restrict our consideration to resultantal varieties that
correspond to essential collections. Hence, in what follows, we
only consider essential collections $\Sigma_1,\ldots,\Sigma_I$,
such that $\Sigma_1+\ldots+\Sigma_I$ is not contained in a
hyperplane.
\begin{lemma} \label{essent1} Let $\Sigma_i\subset \Z^N,
i=1,\ldots,I$, be finite sets.
\newline 1) There exists the minimal subset
$\{i_1,\ldots,i_J\}\subset\{1,\ldots,I\}$, such that
$$\codim(\Sigma_{i_1},\ldots,\Sigma_{i_J})=\max\limits_{\{j_1,\ldots,j_p\}\subset\{1,\ldots,I\}}\codim(\Sigma_{j_1},\ldots,\Sigma_{j_p}).$$
In particular, the collection $\Sigma_{i_1},\ldots,\Sigma_{i_J}$
is essential.
\newline 2) $R(\Sigma_1,\ldots,\Sigma_I)=
p^{(-1)}\bigl(R(\Sigma_{i_1},\ldots,\Sigma_{i_J})\bigr)$, where
$p:\C[\Sigma_1]\oplus\ldots\oplus\C[\Sigma_I]
\to\C[\Sigma_{i_1}]\oplus\ldots\oplus\C[\Sigma_{i_J}]$ is the
natural projection.
\end{lemma}
\textsc{Proof} of Part 1 for $I=N+1$ is given in \cite{sturmf},
Section 1, and can be extended to the general case. Part 2
follows from Part 1 and Lemma \ref{sturm1}. See e.g. \cite{E4} for
details.
\begin{exa}If $\Sigma_1=\Sigma_2$ is a segment in $\R^2$, and $\Sigma_3$ is a polygon,
then $(\Sigma_1,\Sigma_2)$ is the minimal essential subcollection
in the collection $(\Sigma_1,\Sigma_2,\Sigma_3)$, and $\codim
R(\Sigma_1,\Sigma_2,\Sigma_3)=\codim(\Sigma_1,\Sigma_2)=1$.
\end{exa}

\subsection{Cayley trick for intersection numbers}\label{sscay}

In this subsection we study the multiplicity of a 0-dimensional
resultantal singularity in terms of Newton polyhedra. Let
$\Sigma_1,\ldots,\Sigma_I$ be an essential collection of finite
sets in $\Z^N$, such that their sum
$(\Sigma_1+\ldots+\Sigma_I)\times\{1\}$ generates the lattice
$\Z^N\oplus\Z$ (recall that we can impose both of these
restrictions without loss of generality). A point in the space
$\C[\Sigma_1]\oplus\ldots\oplus\C[\Sigma_I]$ is a collection of
Laurent polynomials of the form $(\sum\limits_{a\in\Sigma_1}
y_{(a,1)} t^a, \ldots, \sum\limits_{a\in\Sigma_I} y_{(a,I)}
t^a)$, where $y_m,\; m\in M=\{ (a,i)\;|\;
a\in\Sigma_i,\;i=1,\ldots,I\}$, is the natural system of
coordinates in $\C[\Sigma_1]\oplus\ldots\oplus\C[\Sigma_I]$. A
germ of a holomorphic mapping
$f:(\C^n,0)\to(\C[\Sigma_1]\oplus\ldots\oplus\C[\Sigma_I],0)$ is
given by its components $y_m=f_m,\;m\in M$, in this coordinate
system. We discuss the following problem: how to compute the
intersection number
$$(\mbox{graph of}\; f) \circ \bigl(R(\Sigma_1,\ldots,\Sigma_I)\times\C^n\bigr)\eqno (*)$$ in terms of
the Newton polyhedra $\Delta_{f_m},\;m\in M$, under the assumption
that the leading coefficients of $f_m$ are in general position?

In what follows, we use notation introduced in Subsection
\ref{sstor}. We choose a simple fan $\Gamma$ compatible with the
convex hulls $\Span(\Sigma_1),\ldots,$ $\Span(\Sigma_I)$. The
toric variety $\T^{\Gamma}$ carries the line bundle
$\mathcal{B}_{\Span(\Sigma_i)}$ with the distinguished section
$s_{\Span(\Sigma_i)}$. We pull it back to the product
$\T^{\Gamma}\times\C^n$, and denote the germ of its section
$s_{\Span(\Sigma_i)}\cdot\sum\limits_{a\in\Sigma_i} f_{a,i} t^a$
on the pair $(\T^{\Gamma}\times\C^n, \T^{\Gamma}\times\{0\})$ by
$s_i$. The following theorem allows us to discuss Newton
polyhedra, leading coefficients and topology of the sections
$s_1,\ldots,s_I$ instead of those of the mapping $f$.
\begin{theor}[\cite{E3}] \label{indind}
1) Under the above assumptions, the intersection number $(*)$
equals $\ind ([s_1]\cap\ldots\cap [s_I])$. In particular, these
intersection numbers make sense simultaneously.
\newline 2) The Newton polyhedron of the section ${s_i}$ equals $$\Delta_i=\conv\left(\bigcup\limits_{a\in\Sigma_i}
\{a\}\times\Delta_{f_{(a,i)}}\right)\subset\R^N\times\R^n_+.$$ 3)
Every leading coefficient of $s_i$ is a leading coefficient of
one of the components $f_{(a,i)}$.
\end{theor}
Actually, Part 1 is valid for an arbitrary continuous map $f$,
the assumption of holomorphicity is redundant. Note that Part 1
only makes sense for $n=I-N$ (otherwise, the intersection numbers
are not defined).
\begin{exa}The matrix $A$ from Example \ref{exadet} can be considered as a germ of a mapping
$f:\C^2\to\C[S]\oplus\C[S]\oplus\C[S]$, where $S=\{(0,0),\,
(0,1),\, (1,0)\}$, see Lemma \ref{detviarez}. In this case we have
$\Delta_1=\Delta_2=\Delta_3=\bigl([0,1]\times\R^2_+\bigr)\setminus\Sigma$,
where $\Sigma$ is shown on the picture of Example \ref{exadet}.
Part 1 states that the multiplicity of degeneration of $A$ at the
origin is equal to the intersection number of the divisors
$\lambda a_{1,i}+\mu a_{2,i}=0$ in $\CP^1\times\C^2$, where
$\lambda:\mu$ are the standard coordinates on $\CP^1$, and $i$
runs over $1,2,3$.
\end{exa}
In particular, computing the intersection number $\ind
([s_1]\cap\ldots\cap [s_I])$ on the toric variety
$\T^{\Gamma}\times\C^n$ by Theorem \ref{relbernst}, we have the
following corollary for the mapping $f$.
\begin{sledst}[\cite{E2}, \cite{E3}]\label{corolres1} The intersection number $(*)$ equals $I!$ times the mixed
volume of pairs
$\Bigl(\conv(\Sigma_i)\times\R^n_+,\Delta_i\Bigr)$, provided that
the leading coefficients of the components $f_{(a,i)}$ are in
general position, and the difference of polyhedra in each of
these pairs is bounded.
\end{sledst}

\textsc{Proof of Theorem \ref{indind}.} Parts 2 and 3 follow from
definitions.

Proof of Part 1 (we use notation introduced in Subsection
\ref{sskb}): the line bundle $\mathcal{B}_{\Span(\Sigma_i)}$
lifted from the toric variety $\T^{\Gamma}$ to the product
$M=\T^{\Gamma}\times(\C[\Sigma_1]\oplus\ldots\oplus\C[\Sigma_I])$
admits \textit{the tautological section} $\tilde s_i$, such that
$\tilde
s_i|_{\T^{\Gamma}\times\{(l_1,\ldots,l_I)\}}=s_{\Span(\Sigma_i)}\cdot
l_i$ for every point
$(l_1,\ldots,l_I)\in\C[\Sigma_1]\oplus\ldots\oplus\C[\Sigma_I]$.
The cycle $[R(\Sigma_1,\ldots,\Sigma_I)]$ is the image of the
intersection $[\tilde s_1]\cap\ldots\cap[\tilde s_I]$ under the
projection $\pi$ of the product $M$ to the second factor, the
section $s_i$ is the inverse image of $\tilde s_i$ under the map
$(\id,f):\T^{\Gamma}\times\C^n\to M$, and the diagram
$$\begin{matrix}
   \T^{\Gamma}\times\C^n & \stackrel{(\id,f)}{\longrightarrow} & M \\
  \downarrow &   & \phantom{\pi}\downarrow\pi \\
  \C^n & \stackrel{f}{\longrightarrow} & \C[\Sigma_1]\oplus\ldots\oplus\C[\Sigma_I]
\end{matrix}$$
shows that $\ind (\id,f)^*\bigl([\tilde s_1]\cap\ldots\cap[\tilde
s_I]\bigr)=\ind\pi_*\bigl([\tilde s_1]\cap\ldots\cap[\tilde
s_I]\bigr)\cap f_*[\C^n].$ Indeed,
$$
\pi_* ([\tilde s_1]\cap \ldots \cap [\tilde s_{I}]) \cap
f_*[\C^n] =\pi_* \Bigl([\tilde s_1]\cap \ldots \cap [\tilde s_{I}]
\cap \pi^* f_* [\C^n]\Bigr) =
$$
$$
=\pi_* \Bigl([\tilde s_1]\cap \ldots \cap [\tilde s_{I}] \cap
(\id, f)_* [\T^{\Gamma} \times \C^n]\Bigr) = \pi_*(\id, f)_*(\id,
f)^* ([\tilde s_1]\cap \ldots \cap [\tilde s_{I}]).\; \Box
$$

\subsection{Toric resolutions of resultantal singularities}\label{sstorrez}
We construct a toric resolution for a
positive-dimensional resultantal singularity
$f^{(-1)}\bigl(R(\Sigma_1,\ldots,\Sigma_I)\bigr)$, such that the
leading coefficients of the components $f_{a,i}$ of the mapping
$f$ are in general position.

The faces of the positive orthant $\R^n_+$ form a fan that we
denote by $P$. Adopting notation of Theorem \ref{indind}, we
choose a simple subdivision $\Gamma'$ of the fan $\Gamma\times P$,
compatible with the Newton polyhedra $\Delta_{i},\,
i=1,\ldots,I$. The function $\sum\limits_{a\in\Sigma_i} f_{a,i}
t^a$ on the torus $\CC^N\times\CC^n$, that can be regarded as the
maximal torus of the variety $\T^{\Gamma'}$, defines a holomorphic
section $r_i=s_{\Delta_{i}}\cdot\sum\limits_{a\in\Sigma_i}
f_{a,i} t^a$ of the line bundle $\mathcal{B}_{\Delta_{i}}$ on
this variety.
\begin{defin}[\cite{E4}] The set $\{r_1=\ldots=r_I=0\}\subset\T^{\Gamma'}$, together with its projection $\pi$
to $\C^n$, is a \textit{toric resolution} of the singularity
$f^{(-1)}R(\Sigma_1,\ldots,\Sigma_I)$.\end{defin} This definition
is motivated by the following lemma.
\begin{lemma}[\cite{E4}] \label{lres1} Let $\Sigma_1,\ldots,\Sigma_I$ be an
essential collection of finite sets in $\Z^N$, such that their sum
$(\Sigma_1+\ldots+\Sigma_I)\times\{1\}$ generates the lattice
$\Z^N\oplus\Z$. \newline 1) We have
$\pi\Bigl(\{r_1=\ldots=r_I=0\}\Bigr)=f^{(-1)}R(\Sigma_1,\ldots,\Sigma_I)$.
\newline 2) If the Newton polyhedra of the components $f_{a,i}$ of the map $f$
intersect all coordinate axes in $\R^n$, and the leading
coefficients of the components are in general position, then
$\{r_1=\ldots=r_I=0\}$ is smooth, and the topological degree of
the map $\pi:\{r_1=\ldots=r_I=0\}\to
f^{(-1)}R(\Sigma_1,\ldots,\Sigma_I)$ equals 1. In particular,
$f^{(-1)}R(\Sigma_1,\ldots,\Sigma_I)$ is a resultantal
singularity.
\end{lemma}
The proof is the same as for toric resolutions of complete
intersections. A sufficient condition of general position is
that, for every collection of weights, assigned to the variables
$x_1,\ldots,x_n$ and $t_1,\ldots,t_n$, such that the weights of
$x_1,\ldots,x_n$ are positive, the polynomial equations
$(\sum\limits_{a\in\Sigma_i} f_{a,i} t^a)^l=0,\; i=1,\ldots,I$,
define a non-degenerate subvariety in the torus
$\CC^N\times\CC^n$. Let $\mathcal{Q}$ be 1 if the sets
$\Sigma_1,\ldots,\Sigma_I$ are contained in the standard simplex,
and let it be 0 otherwise.
\begin{theor}[\cite{E4}] \label{smoothres} If $n\leqslant 2(I-N)+\mathcal{Q}$, the convex hulls of $\Sigma_1,\ldots,\Sigma_I$
have the same dual fan, and this fan is simple, then, under the
assumptions of Lemma \ref{lres1}(2), the map
$\pi:\{r_1=\ldots=r_I=0\}\to f^{(-1)}R(\Sigma_1,\ldots,\Sigma_I)$
is a diffeomorphism outside the origin $0\in\C^n$.
\end{theor}
In particular, $f^{(-1)}R(\Sigma_1,\ldots,\Sigma_I)$ is an
isolated resultantal singularity in this case. See \cite{E4} for a
more general statement with no assumptions on the dual fans of
convex hulls of $\Sigma_i$.

\textsc{Proof.} Let $S$ be the set of all collections
$(\varphi_1,\ldots,\varphi_I),\; \varphi_i\in\C[\Sigma_i]$, such
that the sections $s_{\conv(\Sigma_i)}\cdot\varphi_i$ of the line
bundles $\mathcal{I}_{\conv(\Sigma_i)},\; i=1,\ldots,I$, have at
least two different common zeros or a multiple common zero in
$\T^{\Gamma}$ (recall that "$\conv$" stands for the convex hull).
If $f^{(-1)}(S)=\{0\}$ and the set $\{r_1=\ldots=r_I=0\}$ is
smooth, then the map $\pi:\{r_1=\ldots=r_I=0\}\to
f^{(-1)}R(\Sigma_1,\ldots,\Sigma_I)$ is a diffeomorphism outside
the origin $0\in\C^n$. Hence, Theorem \ref{smoothres} is a
corollary of the following two lemmas. $\Box$

\begin{lemma}[\cite{E4}] \label{smoothres1} Under the assumptions of Theorem
\ref{smoothres}, we have $$\codim S\geqslant 2(I-N)+\mathcal{Q}.$$
\end{lemma}
Since a common root of the sections
$s_{\conv(\Sigma_i)}\cdot\varphi_i,\, i=1,\ldots,I$, appears in
codimension $I-N$, it is natural to expect two common roots in
codimension at least $2(I-N)$. If, in addition, the sets
$\Sigma_1,\ldots,\Sigma_I$ are contained in the standard simplex,
then the sections $s_{\conv(\Sigma_i)}\cdot\varphi_i,\,
i=1,\ldots,I$, vanish at every point of the line through the two
common roots, which increases codimension by 1. See \cite{E4} for
a formal proof.
\begin{lemma}[\cite{E4}, \cite{E5}] If the Newton polyhedra of the components of the map $f$
intersect all coordinate axes in $\R^n$, and the leading
coefficients of the components of $f$ satisfy a certain condition
of general position, then the image of $f$ intersects the variety
$S$ properly outside of the origin.
\end{lemma}
In particular, if $n\leqslant\codim S$, then $f^{(-1)}(S)=\{0\}$.
See \cite{E4} or \cite{E5} for the proof. Note that Theorem
\ref{smoothres} and Lemma \ref{smoothres1} provide a sharp
estimate for the codimension of the singular locus of resultantal
singularities and varieties, under the following assumption.
\begin{conjec} If the mixed volume of an essential collection of
integer polyhedra $\Sigma_1,\ldots,\Sigma_m$ in $\R^m$ is equal to
$1/m!$ (i.e. the minimal possible one), then, for a certain
automorphism $A:\Z^m\to\Z^m$ and vectors $a_1,\ldots,a_m$ in
$\Z^m$, the polyhedra $A(\Sigma_1)+a_1,\ldots,A(\Sigma_m)+a_m$
are contained in the standard simplex.
\end{conjec}
To the best of my knowledge, this fact is proved (by G. Gusev)
only under an additional assumption $\dim\Sigma_1=m$ so far.

\subsection{Proofs of results of Section 1}\label{ss3proof}

\textsc{Proof of Theorem \ref{thsturmf}.} If $I=N+1$ and the
collection $\Sigma_1,\ldots,\Sigma_I\subset\Z^N$ is essential,
then we have $R(\Sigma_1,\ldots,\Sigma_I)=\{\Rez=0\}$. Let
$f:\C\to\C[\Sigma_1]\oplus\ldots\oplus\C[\Sigma_I]$ be a germ of
a holomorphic mapping with components $f_m(t)=c_m
t^{\gamma_m}+\ldots$, where $c_m$ are generic non-zero complex
numbers, $m\in M,\; M=\{(a,i)|a\in \Sigma_i\}$. Then the
intersection number \ref{sscay}$.(*)$ equals
$\Delta_{\Rez}(\gamma),$ where $\Delta_{\Rez}(\cdot)$ is the
support function of the Newton polyhedron $\Delta_{\Rez}$, and
$\gamma\in(\R^M)^*$ is the covector with coordinates $\gamma_m,\;
m\in M$. We can compute it by Corollary \ref{corolres1}. $\Box$

\textsc{Proof of Theorem \ref{matrzeta}.} By Lemma
\ref{detviarez}, determinantal singularities can be regarded as
resultantal singularities. The construction of the toric
resolution for resultantal singularities (Theorem
\ref{smoothres}) implies Part 1, and reduces Part 2 to Theorem
\ref{relbkh}. Part 3 can be reduced to Part 2 by means of this
construction, in the same way as it is done in \cite{E1} for
complete intersections by means of their toric resolutions (the
idea is to deform the differential form $\omega$ to the
differential of a function, preserving its radial index and
Newton polyhedron). $\Box$

\textsc{Proof of Theorem \ref{egvol}.} Lemma \ref{detviarez}
implies that the multiplicity of the collection of matrices is a
special case of the multiplicity of a 0-dimensional resultantal
singularity, treated in Corollary \ref{corolres1}. $\Box$

\textsc{Proof of Theorem \ref{finoka1}.} The proof is not exactly
the same as for Theorem \ref{egvol}, because the straightforward
application of Corollary \ref{corolres1} would give an answer
involving the Newton polyhedra of the partial derivatives of the
germs $f_j$.

We apply induction on $n$. The desired index equals the
multiplicity of the collection
of matrices $$\begin{pmatrix}x_1w_1 & \ldots & x_nw_n \\
x_1\frac{\partial f_1}{\partial x_1} & \ldots & x_n\frac{\partial
f_1}{\partial x_n} \\ \hdotsfor{3} \\ x_1\frac{\partial
f_k}{\partial x_1} & \ldots & x_n\frac{\partial f_k}{\partial x_n}
\end{pmatrix},\bigl(f_1,\ldots,f_k\bigr)$$
minus the sum of the indices of the 1-form $w_1 dx_1 + \ldots +
w_n dx_n$ restricted to the complete intersections
$f_1=\ldots=f_k=x_{i_1}=\ldots=x_{i_m}=0$, over all non-empty
subsets $\{i_1,\ldots,i_m\}\subset\{1,\ldots,n\}$. The latter
indices can be computed in terms of Newton polyhedra by
induction, and the multiplicity of the collection of matrices
equals
the multiplicity of the collection $$\begin{pmatrix} x_1w_1 & \ldots & x_nw_n \\
x_1\frac{\partial f_1}{\partial x_1}+a_{1,1}f_1 & \ldots &
x_n\frac{\partial f_1}{\partial x_n}+a_{1,n}f_1
 \\ \hdotsfor{3} \\ x_1\frac{\partial f_k}{\partial
x_1}+a_{k,1}f_k & \ldots & x_n\frac{\partial f_k}{\partial
x_n}+a_{k,n}f_k\end{pmatrix},\bigl(f_1,\ldots,f_k\bigr),$$ where
$a_{i,j}$ are generic complex coefficients. Lemma \ref{detviarez}
implies that this multiplicity is a special case of the
multiplicity of a 0-dimensional resultantal singularity. Theorem
\ref{indind} represents it as the intersection number of certain
divisors on a toric variety, which can be computed/estimated in
terms of Newton polyhedra by Corollary \ref{volveryconv}. The
resulting answer does not involve the Newton polyhedra of the
partial derivatives of $f_i$ because the Newton polyhedron of
$x_j\frac{\partial f_i}{\partial x_j}+a_{i,j}f_i$ equals
$\Delta_{f_i}$ (in contrast to $\Delta_{x_j\frac{\partial
f_i}{\partial x_j}}$) $.\; \Box$

\subsection{Proof of Theorem \ref{thint}}\label{ss4}

This proof is purely combinatorial, it does not rely upon
material of Sections \ref{ss2} and \ref{ss3}.

\textbf{Lattice points of sums of polyhedra.} We prove the
equality
$$(A\cap\Z^q)+(B\cap\Z^q)=(A+B)\cap\Z^q$$ for a certain class of
bounded integer polyhedra $A,B$ in $\R^q$ (cf. \cite{oda1}).
\begin{defin}
A collection of rational cones $C_1,\ldots,C_p$ in $\R^q$ is said
to be $\Z$-\textit{transversal}, if $\sum\dim C_i=q$, and the set
$\Z^q\cap\bigcup_i {C_i}$ generates the lattice $\Z^q$.
\end{defin}
\begin{defin}
A collection of fans $\Phi_1,\ldots,\Phi_p$ in $\R^q$ is said to
be $\Z$-\textit{transversal with respect to shifts} $c_1\in\R^q,
\ldots, c_p\in\R^q$, if every collection of cones $C_1\in\Phi_1,
\ldots, C_p\in\Phi_p$, such that the intersection
$(C_1+c_1)\cap\ldots\cap (C_p+c_p)$ consists of one point, is
$\Z$-transversal.
\end{defin}
\begin{theor} \label{transvpoly}
If the dual fans of bounded integer polyhedra $A_1,\ldots,A_p$ in
$\R^q$ are $\Z$-transversal with respect to certain shifts
$c_1\in(\R^q)^*, \ldots, c_p\in(\R^q)^*$, and
$\dim(A_1+\ldots+A_p)=q$, then
$$(A_1\cap\Z^q)+\ldots+(A_p\cap\Z^q)=(A_1+\ldots+A_p)\cap\Z^q.$$
\end{theor}

\textsc{Proof.} Consider covectors $c_1\in(\R^q)^*, \ldots,
c_p\in(\R^q)^*$ as linear functions on the polyhedra
$A_1\subset\R^q,\ldots,A_p\subset\R^q$ respectively, and denote
their graphs in $\R^q\oplus\R^1$ by $\Gamma_1,\ldots,\Gamma_p$.
Denote the projection $\R^q\oplus\R^1\to\R^q$ by $\pi$, and
denote the ray $\{ (0,\ldots,0,t)\; |\;
t<0\}\subset\R^q\oplus\R^1$ by $L_-$.

Each bounded $q$-dimensional face $B$ of the sum
$\Gamma_1+\ldots+\Gamma_p+L_-$ is the sum of certain faces
$B_1,\ldots,B_p$ of polyhedra $\Gamma_1+L_-,\ldots,\Gamma_p+L_-$,
and $\Z$-transversality with respect to shifts $c_1\in(\R^q)^*,
\ldots, c_p\in(\R^q)^*$ implies that
$$\bigl(\pi(B_1)\cap\Z^q\bigr)+\ldots+\bigl(\pi(B_p)\cap\Z^q\bigr)=\pi(B_1+\ldots+B_p)\cap\Z^q.$$
Since the projections of bounded $q$-dimensional faces of the sum
$\Gamma_1+\ldots+\Gamma_p+L_-$ cover the sum $A_1+\ldots+A_p$, it
satisfies the same equality:
$$(A_1\cap\Z^q)+\ldots+(A_p\cap\Z^q)=(A_1+\ldots+A_p)\cap\Z^q.\;\Box$$

\begin{sledst} \label{minimax1}
Let $S\subset\R^q$ be the standard $q$-dimensional simplex, let
$l_1,\ldots,l_p$ be linear functions on $S$ with graphs
$\Gamma_1,\ldots,\Gamma_p$, and let $l$ be the maximal
piecewise-linear function on $pS$, such that its graph $\Gamma$
is contained in the sum $\Gamma_1+\ldots+\Gamma_p$. Then, for
each integer lattice point $a\in pS$, the value $l(a)$ equals the
maximum of sums $l_1(c_1)+\ldots+l_p(c_p)$ over all $p$-tuples
$(c_1,\ldots,c_p)$ of vertices of $S$, such that
$c_1+\ldots+c_p=a$.
\end{sledst}
\textsc{Proof.} Denote the projection $\R^q\oplus\R^1\to\R^q$ by
$\pi$. A $q$-dimensional face $B$ of $\Gamma$, which contains the
point $\bigl(a,l(a)\bigr)\in\R^q\oplus\R^1$, can be represented
as a sum of faces $B_i$ of simplices $\Gamma_i$. Since
$\pi(B_1),\ldots,\pi(B_p)$ are faces of the standard simplex,
their dual fans are $\Z$-transversal with respect to a generic
collection of shifts, and, by Theorem \ref{transvpoly},
$$\bigl(\pi(B_1)\cap\Z^q\bigr)+\ldots+\bigl(\pi(B_p)\cap\Z^q\bigr)=\pi(B)\cap\Z^q.$$
In particular, $a=c_1+\ldots+c_p$ for some integer lattice points
$c_i\in\pi(B_i)$, which implies $l(a)=l_1(c_1)+\ldots+l_p(c_p)$.
$\Box$

\begin{rem}
In particular, if the functions $l_1,\ldots,l_p$ are in general
position, then all $C_{p+q}^{q}$ integer lattice points in the
simplex $pS$ are projections of vertices of $\Gamma$. In the
tropical language, this is a well-known fact that $p$ generic
tropical hyperplanes in the space $\R^q$ subdivide it into
$C_{p+q}^{q}$ pieces.
\end{rem}

\begin{exa}
If $S$ in the formulation of Corollary \ref{minimax1} is not the
standard simplex, then the statement is not always true. For
example, consider $$S=\bigl\{|x|+|y|\leqslant 1\bigr\},\;
l_1(x,y)=x+y,\; l_2(x,y)=x-y,\; a=(1,0).$$

More generally, if $S_1,\ldots,S_p$ is an essential collection of
polyhedra in $\R^q,\, q>1$, that cannot be represented as a
collection of shifted faces of an elementary integer simplex, then
there exist concave piecewise linear functions $l_j:S_j\to\R$
with integer domains of linearity, such that the value of the
corresponding function $l(a)=\max\limits_{a_j\in S_j,\;
a_1+\ldots+a_p=a}\sum l_j(a_j)$ on $S_1+\ldots+S_p$ at some
integer point $a$ is strictly greater than $\max\limits_{a_j\in
S_j\cap\Z^q\atop a_1+\ldots+a_p=a}\sum l_j(a_j)$. That is why we
cannot extend Theorem \ref{thint} to mixed volumes of Newton
polyhedra, related to resultantal singularities, other than
determinantal ones.
\end{exa}

\textbf{Proof of Theorem \ref{thint}.} The desired statement
follows from Lemmas \ref{volint} and \ref{intsum} below.

\begin{lemma} \label{volint}
For pairs of integer polyhedra $A_i\in\mathcal{M}_{\Gamma}$, we
have
$$n!\Vol_{\Gamma}(A_1,\ldots,A_n)=\sum_{1\leqslant i_1<\ldots<i_p\leqslant n}
(-1)^{n-p} I(A_{i_1}+\ldots+A_{i_p})+(-1)^n.$$
\end{lemma}
The proof is the same as for the classical mixed volume (see e.g.
\cite{pkhovvol}).
\begin{lemma} \label{intsum}
For pairs of polyhedra $B_{i,j}\in\mathcal{M}_{\Gamma}$,
$i=1,\ldots,n$, $j=1,\ldots,p$, we have
$$I\bigl(B_{1,1}*\ldots*B_{n,1}+\ldots+B_{1,p}*\ldots*B_{n,p}\bigr)=$$
$$=\sum_{a_1+\ldots+a_n=p\atop a_1\geqslant 0,\ldots, a_n\geqslant 0}I\Bigl(
\bigvee_{J_1\sqcup\ldots\sqcup J_n=\{1,\ldots,p\}\atop
|J_1|=a_1,\ldots,|J_n|=a_n} \sum_{i=1,\ldots,n\atop j\in J_i}
B_{i,j} \Bigr).$$
\end{lemma}
\textsc{Proof.} Every integer lattice point that participates in
the left hand side, is contained in the plane
$\{(a_1,\ldots,a_{n-1})\}\times\R^m\subset\R^{n-1}\oplus\R^m$ for
certain non-negative integer numbers $a_1,\ldots,a_n$, which sum
up to $p$. Thus, it is enough to describe the intersection of the
pair
$\bigl(B_{1,1}*\ldots*B_{n,1}+\ldots+B_{1,p}*\ldots*B_{n,p}\bigr)$
with each of these planes, using the following fact. $\Box$

\begin{lemma} Suppose that polyhedra $\Delta_{i,j}\subset\R^m$ have the same support cone for $i=1,\ldots,n$,
$j=1,\ldots,p$. Then, for each $n$-tuple of non-negative integer
numbers $a_1,\ldots,a_n$ which sum up to $p$,
$$\Bigl(\{(a_1,\ldots,a_{n-1})\}\times\R^m\Bigr)\, \cap\,
\bigl(
\Delta_{1,1}*\ldots*\Delta_{n,1}+\ldots+\Delta_{1,p}*\ldots*\Delta_{n,p}\bigr)=$$
$$=\{(a_1,\ldots,a_{n-1})\}\times\Bigl(\bigvee_{J_1\sqcup\ldots\sqcup J_n=\{1,\ldots,p\}\atop |J_1|=a_1,\ldots,|J_n|=a_n}
\sum_{i=1,\ldots,n\atop j\in J_i}
\Delta_{i,j}\Bigr)\subset\R^{n-1}\oplus\R^m.$$
\end{lemma}

\textsc{Proof.} For every hyperplane $L\subset\R^m$, denote the
projection $\R^{n-1}\oplus\R^m\to\R^{n-1}\oplus\R\, $ along $\,
\{0\}\oplus L$ by $\pi_L$. It is enough to prove that the images
of the left hand side and the right hand side under $\pi_L$
coincide for every $L$. To prove it, apply Corollary
\ref{minimax1}, assuming that $q=n-1,\; a=(a_1,\ldots,a_{n-1})$,
and $\Gamma_j$ is the maximal bounded face of the projection
$\pi_L\bigl(\Delta_{1,j}*\ldots*\Delta_{n,j}\bigr)$ for every
$j=1,\ldots,p$. $\Box$


\end{document}